\documentstyle{amsppt}
\magnification=\magstep1
\hsize = 150truemm
\vsize = 230truemm
\hcorrection{4.8truemm}
\vcorrection{8truemm}
\TagsOnRight
\NoBlackBoxes

\define\fine{\operatorname{fine}}
\redefine\ge{\geqslant}
\redefine\le{\leqslant}

\define\QB{\operatorname{QB}}

\define\wh{\widehat}

\define\CC{{\Bbb C}}

\redefine\epsilon{\varepsilon}
\redefine\phi{\varphi}
\redefine\Re{\operatorname{Re}}
\redefine\Im{\operatorname{Im}}

\document

\medskip
\centerline{\bf PLURISUBHARMONIC AND HOLOMORPHIC FUNCTIONS}
\centerline{\bf RELATIVE TO THE PLURIFINE TOPOLOGY}

\medskip
\medskip

\centerline{\smc Mohamed El Kadiri, Bent Fuglede, and Jan Wiegerinck}
\topmatter

\medskip\medskip
 \abstract A weak and a strong concept of plurifinely plurisubharmonic
and plurifinely holomorphic functions are introduced. Strong will
imply weak. The weak concept is studied further. A function $f$ is
weakly plurifinely plurisubharmonic if and only if $f\circ h$ is
finely subharmonic for all complex affine-linear maps $h$. As a
consequence, the regularization in the plurifine topology of a
pointwise supremum of such functions is weakly plurifinely
plurisubharmonic, and it differs from the pointwise supremum at most
on a pluripolar set.  Weak plurifine plurisubharmonicity and weak
plurifine holomorphy are preserved under composition with weakly
plurifinely holomorphic maps.
 \endabstract

\footnote""{2000 {\it Mathematics Subject Classification.} 31C40, 32U15}

\footnote""{{\it Keywords and phrases.} plurifinely plurisubharmonic
function, finely subharmonic function, plurifine topology}
\endtopmatter

\bigskip\bigskip
\heading{\bf 1. Introduction}\endheading

\smallskip The {\it pluri\-fine topology} $\Cal F$ on $\Bbb C^n$ was
briefly introduced in \cite{F6} as the weakest topology in which all
pluri\-sub\-har\-monic functions are con\-tinu\-ous, in analogy with
the H.\ Cartan fine topology on $\Bbb R^n$, in particular on $\Bbb
C\cong\Bbb R^2$. For comments on this choice of \lq\lq fine'' topology on
$\Bbb C^n$, see \cite{F6}. The plurifine topology $\Cal F$ is
clearly biholo\-morph\-ic\-ally invariant. Furthermore, $\Cal F$ is
locally connected, as shown in \cite{EW1}, \cite{EW2}, where also
further properties of \,$\Cal F$ are given. Much as in \cite{ElKa}, \cite{EW2},
\cite{EW3} we begin by considering (in Definition 2.2, resp.\ 2.6) two
concepts of pluri\-finely pluri\-sub\-har\-monic (resp.\ pluri\-finely
ho\-lo\-morph\-ic) functions---a strong concept defined by $\Cal
F$-local uniform approximation with pluri\-sub\-har\-monic (resp.\
ho\-lo\-morph\-ic) functions, and a weak concept defined by
restriction to complex lines. We thereby draw on the theory of {\it
finely sub- or super\-harmonic} and {\it finely ho\-lo\-morph\-ic}
functions defined on finely open subsets of \,$\Bbb C$, cf.\
\cite{F1}, \cite{F3}, \cite{F7}. The pluri\-fine topology $\Cal F$ on $\Bbb
C^n$ induces on each complex line $L$ in $\Bbb C^n$ the Cartan fine
topology on $L\cong\Bbb C$ (Lemma 2.1). In analogy with ordinary
pluri\-sub\-har\-monic functions, the weakly $\Cal
F$-pluri\-sub\-har\-monic functions $f$ may be characterized by being
$\Cal F$-upper semicontinuous and such that $f\circ h$ is $\Bbb
R^{2n}$-finely subharmonic (or identically $-\infty$ in some fine
component of its domain of definition) for every $\Bbb C$-affine-linear
bijection $h$ of $\Bbb C^n$ (Theorem 3.1).

The concepts of strongly $\Cal F$-pluri\-sub\-har\-monic and strongly
$\Cal F$-holomorphic functions on an $\Cal F$-domain
$\Omega\subset\Bbb C^n$ are obviously biholo\-morph\-ic\-ally
invariant. We show that the same holds for the weak concepts (Theorem
4.6), cf.\ \cite{EW3}. We do not know whether the strong and the weak
concepts are actually the same. The weak concepts are closed under
$\Cal F$-locally uniform convergence, and seem altogether to be more
useful, cf.\ \cite{EW2}.

The convex cone of all weakly $\Cal F$-pluri\-sub\-har\-monic
functions on $\Omega$ is stable under pointwise infimum for lower
directed families and under pointwise supremum for finite
families. The above characterization of weakly $\Cal
F$-pluri\-sub\-har\-monic functions allows us to answer questions
posed by the first named author in \cite{ElKa}.  Namely, for any $\Cal
F$-locally upper bounded family of \,weakly $\Cal
F$-pluri\-sub\-har\-monic functions $f_\alpha$ on $\Omega$, the $\Cal
F$-upper semi\-con\-tinu\-ous regularization $f^*$ of the pointwise
supremum $f=\sup_\alpha f_\alpha$ is likewise weakly $\Cal
F$-pluri\-sub\-har\-monic (Theorem 3.9), and the exceptional set
$\{f<f^*\} $ is pluripolar, as expected from a theorem of Bedford and
Taylor \cite{BT1, Theorem 7.1}. Furthermore, there is a removable
singularity theorem for weakly $\Cal F$-pluri\-sub\-har\-monic
functions (Theorem 3.7), and likewise for $\Cal F$-holomorphic functions (Corollary 3.8).

In the final Section 4 we show that the concepts of weakly $\Cal
F$-pluri\-sub\-har\-monic map and weakly $\Cal F$-holo\-morph\-ic map
are biholo\-morph\-ically invariant, even in a plurifine sense. In
fact, composition with weakly $\Cal F$-holo\-morph\-ic maps preserves
weak $\Cal F$-pluri\-sub\-har\-moni\-city and weak $\Cal F$-holomorphy
(Theorem 4.6).
  
\heading{\bf 2. Definitions and first properties of strongly and
weakly\\ $\Cal F$-pluri\-sub\-har\-monic and $\Cal F$-holo\-morph\-ic
functions}\endheading

The $\Cal F$-interior (pluri\-fine interior) of a set $K\subset\Bbb
C^n$, $n\in\Bbb N$, is denoted by $K'$. It is known that every $\Cal
F$-neighborhood of a point of \,$\Bbb C^n$ contains an $\Cal
F$-neighborhood which is compact in the Euclidean topology---an easy
consequence of \cite{EW2, Theorem 2.3}, pluri\-sub\-har\-monic
functions being upper semi\-con\-tinu\-ous. Henceforth, topological
properties not explicitly referring to the plurifine topology $\Cal F$
or the Cartan fine topology are tacitly understood to refer to the
Euclidean topology. Generalizing known properties of the fine
topology, cf.\ \cite{F8}, we have

\proclaim{Lemma 2.1} {\rm(a)} The plurifine topology $\Cal F$ on $\Bbb
C^n$ induces on every $\Bbb C$-linear subspace $L\cong\Bbb C^k$ of
\,$\Bbb C^n$ the plurifine topology on $L$. Explicitly, for any $\Cal
F$-open set $\Omega\subset\Bbb C^n$ the intersection $L\cap\Omega$ is
$\Cal F$-open in $L$, and so is the orthogonal projection of \,$\Omega$
on $L$.

{\rm(b)} A set $\omega\subset\Bbb C^k$ is $\Cal F$-open in $\Bbb C^k$
if and only if \,$\,\omega\times\Bbb C^{n-k}$ is $\Cal F$-open in $\Bbb
C^n$.
 \endproclaim

\demo{Proof}  For $z=(z_1,\dots,z_n)\in\Bbb C^n$ write
 $$
z'=(z_1,\dots,z_k),\quad z''=(z_{k+1},\dots,z_n).
 $$

 For (a) it suffices to consider the particular subspace
$L_0=\{(z',0''):z'\in\Bbb C^k\}$ which we identify
with $\Bbb C^k$. For any $\Cal F$-open set $\Omega\subset\Bbb C^n$
denote by $\omega$ the part of $\Omega$ in $L_0$.  Consider a point
$a'\in\omega$. According to \cite{EW1, Theorem 2.3} there exists a
pluri\-sub\-har\-monic function $\psi$ on $\Bbb C^n\cong\Bbb
C^k\times\Bbb C^{n-k}$ and neighborhoods $U'$ of $a'$ in $\Bbb C^k$
and $U''$ of $0''$ in $\Bbb C^{n-k}$ such that
 $$
(a',0'')\in\{(z',z'')\in U'\times U'':
\psi(z',z'')>0\}\subset\Omega.\tag2.1
 $$
Define $\phi:\Bbb C^k\to[-\infty,+\infty[\,$ by $\phi(z')=\psi(z',0'')$;
then $\phi$ is pluri\-sub\-har\-monic and
 $$
a'\in\{z'\in U':\phi(z')>0\}\subset\omega.\tag2.2
 $$
Thus $\omega$ is indeed an $\Cal F$-neighborhood of $a'$ in $\Bbb
C^k$.

For each $t\in\Bbb C^{n-k}$ the translate $\Omega_t=\Omega-(0',t)$ of
$\Omega$ is $\Cal F$-open in $\Bbb C^n$. It follows that $\Omega_t\cap
L_0$ is $\Cal F$-open in $L_0$, and so is therefore the union of the
$\Omega_t\cap L_0$, that is, the projection of $\Omega$ on $L_0$.

For (b) we have just shown, in particular, that if
$\Omega:=\omega\times\Bbb C^{n-k}$ is $\Cal F$-open in $\Bbb C^n$ then
$\omega$ is $\Cal F$-open in $\Bbb C^k$. To establish the converse,
suppose that $\omega$ is $\Cal F$-open in $\Bbb C^k$ and let us prove
that every point $a=(a',a'')$ of $\omega\times\Bbb C^{n-k}$ is an
$\Cal F$-inner point of that set. Since $\omega$ is an $\Cal
F$-neighborhood of $a'$ in $\Bbb C^k$ there exists (again by
\cite{EW1, Theorem 2.3}) a pluri\-sub\-har\-monic function $\phi$ on
$\Bbb C^k$ and a neighborhood $U'$ of $a'$ in $\Bbb C^k$ such that
(2.2) holds. The function $\psi$ defined on $\Bbb C^n$ by
$\psi(z',z'')=\phi(z')$ is pluri\-sub\-har\-monic---an easy and
well-known consequence of the definition of
pluri\-sub\-har\-moni\-city \cite{L2, p.\ 306}, or see \cite{K, p.\
62}. Furthermore, (2.1) holds (with $\Omega=\omega\times\Bbb C^{n-k}$
and with $(a',0'')$ replaced by $a$) for any neighborhood $U''$ of
$a''$ in $\Bbb C^{n-k}$. Thus $\omega\times\Bbb C^{n-k}$ is indeed an
$\Cal F$-neighborhood of $a$ in $\Bbb C^n$.
{\hfill$\square$}\enddemo

For a compact set $K\subset\Bbb C^n$ we denote by $S_0(K)$ the convex
cone of all restrictions to $K$ of {\it finite con\-tinu\-ous}
pluri\-sub\-har\-monic functions defined on open subsets of \,$\Bbb
C^n$ containing $K$, and by $S(K)$ the closure of \,$S_0(K)$ in
$C(K,\Bbb R)$ (the continuous functions $K\to\Bbb R$ with the uniform
norm); then $S(K)$ is likewise a convex cone.

\definition{Definition 2.2} Let $\Omega$ denote an $\Cal F$-open
(i.e., pluri\-finely open) subset of \,$\Bbb C^n$.

(i) A function $f:\Omega\to\Bbb R$ is said to be {\it $\Cal F$-cpsh}
if every point of \,$\Omega$ has a compact $\Cal F$-neighborhood $K$
in $\Omega$ such that $f|K\in{S(K)}$.

(ii) A function $f:\Omega\to[-\infty,+\infty[\,$ is said to be {\it
strongly $\Cal F$-pluri\-sub\-har\-monic} if $f$ is the pointwise
limit of a decreasing net of \,$\Cal F$-cpsh
functions on $\Omega$.

(iii) (cf.\ \cite{ElKa, Section 5} \cite{EW2, Definition 5.1}). A function
$f:\Omega\to[-\infty,+\infty[\,$ is said to be {\it weakly $\Cal
F$-pluri\-sub\-har\-monic} if $f$ is $\Cal F$-upper
semi\-con\-tinu\-ous and, for every complex line $L$ in $\Bbb C^n$,
the restriction of $f$ to the finely open subset $L\cap\Omega$ of $L$
is {\it finely hypoharmonic}.
 \enddefinition

See \cite{F1, Definition 8.2 and \S10.4} for finely hypoharmonic
(resp.\ finely sub- or super\-harmonic) functions, and recall that a
function $f$ is finely hypoharmonic on a finely open subset $U$ of
$\Bbb C$ (or of $\Bbb R^N$) if and only if $f$ is finely sub\-harm\-onic
on every fine component of $U$ in which $f\not\equiv-\infty$.  Either
concept strongly or weakly $\Cal F$-pluri\-sub\-har\-monic is an $\Cal
F$-local one (that is, has the sheaf property).

The concept of $\Cal F$-cpsh functions, defined in (i), is an auxiliary
one. Every strongly $\Cal F$-pluri\-sub\-har\-monic function is $\Cal
F$-upper semi\-con\-tinu\-ous (even $\Cal F$-con\-tinu\-ous, see
Theorem 2.4(c) and Proposition 2.5) because every $\Cal F$-cpsh
function is $\Cal F$-conti\-nu\-ous. The class of all strongly, resp.\
weakly, $\Cal F$-pluri\-sub\-har\-monic functions on $\Omega$ is
clearly a convex cone which is stable under pointwise supremum of
finite families. The latter class is furthermore stable under
pointwise infimum for lower directed (possibly infinite) families, and
closed under $\Cal F$-locally uniform convergence in view of \cite{F1,
Lemma 9.6}.  For upper directed families of weakly $\Cal
F$-pluri\-sub\-har\-monic functions, see Theorem 3.9 below.

If $f$ is strongly, resp.\ weakly, $\Cal F$-pluri\-sub\-har\-monic on
$\Omega$ ($\Cal F$-open in $\Bbb C^n$) then the restriction of $f$ to
$L\cap\Omega$ ($L$ a $\Bbb C$-linear subspace $L\cong\Bbb C^k$ of
$\Bbb C^n$) has the same property in $L\cap\Omega$. This follows
easily from Lemma 2.1(a) above.

For $n=1$, $f$ is strongly, resp.\ weakly, $\Cal
F$-pluri\-sub\-har\-monic on $\Omega$ (finely open in $\Bbb C$) if and
only if $f$ is finely hypoharmonic on $\Omega$.  This is obvious in
the weak case. In the strong case, suppose first that $f$ is finite
and finely hypoharmonic on $\Omega$. By the Brelot property \cite{F7,
p.\,248}, every point of $\Omega$ has a compact fine neighborhood $K$
in $\Omega$ such that $f|K\in C(K,\Bbb R)$ ($f$ being finely
continuous by \cite{F1, Theorem 9.10}). Because $f$ is finite and
finely hypoharmonic in the fine interior $K'$ of $K$ we have $f\in
S(K)$ according to \cite{BH, Theorem 4.7}, or see \cite{F5, Theorem
4}, and so $f$ is $\Cal F$-cpsh on $\Omega$. For a general
finely hypoharmonic function $f$ on $\Omega$ write $f=\inf_{n\in\Bbb
N}\max\{f,-n\}$ and note that $\max\{f,-n\}$ is finite and finely
hypoharmonic, cf.\ \cite{F1, Corollary 2, p.\ 84}. Conversely, if $f$
is strongly $\Cal F$-pluri\-sub\-har\-monic we may assume by the same
corollary that $f$ is even $\Cal F$-cpsh. For any compact set
$K\subset\Bbb C$, every function of class $S(K)$ is finite and finely
hypoharmonic on $K'$ according to \cite{F1, Lemma 9.6}. With $K$ as in
(i) this shows that $f$ indeed is finite and finely hypo\-har\-mon\-ic
on $\Omega$.

In the following two theorems we collect some properties of weakly
finely pluri\-sub\-har\-monic functions recently obtained by the
third named author in collaboration with S. El Marzguioui. By an
{\it $\Cal F$-domain} we understand an $\Cal F$-connected $\Cal
F$-open set.

\proclaim{Theorem 2.3} $($\cite{EW2}$)$ Let $f$ be a weakly $\Cal
F$-pluri\-sub\-har\-monic function on an $\Cal F$-domain
$\Omega\subset\Bbb C^n$.

{\rm(a)} If \,$f\not\equiv-\infty$ then $\{z\in\Omega:f(z)=-\infty\}$
has no $\Cal F$-interior point.

{\rm(b)} If \,$f\not\equiv-\infty$ then, for any $\Cal F$-closed set
$E\subset\{z\in\Omega:f(z)=-\infty\}$, $\Cal\Omega\setminus E$ is an
$\Cal F$-domain.

{\rm(c)} If \,$f\le0$ then either $f<0$ or $f\equiv0$.
 \endproclaim

\proclaim{Theorem 2.4} $($\cite{EW3}$)$ Let $f$ be a weakly $\Cal
F$-pluri\-sub\-har\-monic function on an $\Cal F$-open set
$\Omega\subset\Bbb C^n$.
 
{\rm(a)} Every point $z_0\in\Omega$ such that $f(z_0)>-\infty$ has an
$\Cal F$-open $\Cal F$-neighborhood $O\subset\Omega$ on which $f$ can
be represented as the difference $f=\phi_1-\phi_2$ between two bounded
plurisubharmonic functions $\phi_1$ and $\phi_2$ defined on some open
ball $B(z_0,r)$ containing $O$.

{\rm(b)} $r$, $O$, and $\phi_2$ can be chosen independently of \,$f$
provided that $f$ maps $\Omega$ into a prescribed bounded interval
$]a,b[$.

{\rm(c)} $f$ is $\Cal F$-continuous.

{\rm(d)} If \,$\Omega$ is $\Cal F$ connected and $f\not\equiv-\infty$
then $\{z\in\Omega:f(z)=-\infty\}$ is a pluripolar subset of \,$\Bbb C^n$.
 \endproclaim

Assertion (d) amounts to pluripolar sets and weakly $\Cal F$-pluri\-polar
sets (in the obvious sense) being the same. The proofs of (a), (b),
and (c) given below are essentially taken from \cite{EW3}.

\demo{Proof} (a) To begin with, suppose that $f$ is {\it bounded}. We
may then assume that $-1<f<0$, for $f$ maps $\Omega$ into a bounded
interval $\,]a,b[\,$, and hence $\frac{f-b}{b-a}$ maps $\Omega$ into
$\,]-1,0[\,$ and is likewise weakly $\Cal F$-pluri\-sub\-har\-monic.
Let $V \subset \Omega$ be a compact $\Cal{F}$-neighborhood of
$z_0$. Since the complement $\complement V$ of $V$ is pluri-thin at
$z_0$, there exist $ 0<r<1$ and a pluri\-sub\-har\-monic function
$\varphi$ on $B(z_0,r)$ such that
 $$
\limsup_{z \to z_{0},\,z\in \complement V }\varphi(z)<\varphi(z_{0}).
 $$
Without loss of generality we may suppose that $\varphi $ is negative
on $B(z_{0}, r)$ and
$$
\varphi(z)=-1\;\text{on}\;B(z_0,r)\setminus
V\quad\text{and}\;\phi(z_0)=-1/2.
 $$
Hence
 $$
f(z)+\lambda\phi(z)\le-\lambda\qquad\text{for }z\in \Omega\cap
B(z_0,r)\setminus V\text{ and }\lambda>0.\tag2.3
 $$
Now define a function $u_{\lambda}$ on $B(z_0,r)$ by 
 $$
u_\lambda(z)=\cases\max\{-\lambda,f(z)+\lambda\phi(z)\}&\text{for
}z\in \Omega\cap B(z_0,r)\\ -\lambda &\text{for }z\in B(z_0,
r)\setminus V.
 \endcases\tag2.4$$
This definition makes sense because 
$\bigl(\Omega\cap B(z_0,r)\bigr)\bigcup\bigl(B(z_0, r)\setminus
V\bigr)=B(z_0, r)$, and the two definitions agree on $\Omega\cap B(z_{0},
r)\setminus V$ in view of (2.3).

Clearly, $u_{\lambda}$ is weakly $\Cal F$-pluri\-sub\-har\-monic on
$\Omega \cap B(z_0, r)$ and on $B(z_0, r)\setminus V$, hence on all of
$B(z_0, r)$ in view of the sheaf property, cf.\ \cite{EW2}. Since
$u_{\lambda}$ is bounded on $B(z_0, r)$, it follows from \cite{F1,
Theorem 9.8} that $u_{\lambda}$ is subharmonic on each complex line
where it is defined. It is well known that a bounded function, which
is subharmonic on each complex line where it is defined, is
pluri\-sub\-har\-monic, cf. \cite{Le1}, or see \cite{Le2, p.\
24}. Thus, $u_{\lambda}$ is pluri\-sub\-har\-monic on $B(z_0,r)$.

Since $\phi(z_0)=-1/2$, the set $O=\{z\in\Omega:\phi(z)>-3/4\}$ is an
$\Cal{F}$-neighborhood of $z_0$, and because $\phi=-1$ on $B(z_0,
r)\setminus V$ it is clear that $O\subset V\subset\Omega$.

Observe now that $-4\leq f(z)+ 4\phi(z)$ for every $z\in O$. Hence
$f=\phi_1-\phi_2$ on $O$, with $\phi_1=u_4$ and $\phi_2=4\phi$, both
pluri\-sub\-har\-monic on $B(z_0,r)$. Thus $f$ is weakly $\Cal
F$-pluri\-sub\-har\-monic on $O$, which is an $\Cal F$-neighborhood of
$z_0$. It follows that $f$ is $\Cal F$-continuous on $O$ along with
$\phi_1$ and $\phi_2$, provided that $f$ is bounded.  

Without assuming that $f$ be bounded, $f$ remains $\Cal F$-continuous
on $O$ according to (c), proven below. It follows that $f$ is bounded
on some $\Cal F$-neighborhood $U$ of $z_0$ in $\Omega$, and we
therefore have a decomposition of $f$ as required, on some $\Cal
F$-neighborhood (replacing the above $O$) of $z_0$ on
$U\subset\Omega$.

(b) Again we may assume that $-1<f<0$. The set $V$ and the
pluri\-sub\-har\-monic function $\phi$ in the proof of (a) then do not
depend on $f$, and that applies to $\phi_2=4\phi$ as well.

(c) In the remaining case where $f$ may be unbounded (cf.\ the proof
of (a) above), note that $f$ is $\Cal F$-upper semi\-continuous and
$<+\infty$. Choose $c,d\in\Bbb R$ with $d<c$.  Then the set
$\Omega_c=\{z\in\Omega:f(z)<c\}$ is $\Cal F$-open. 
The function $\max\{f,d\}$ is bounded and weakly $\Cal
F$-pluri\-sub\-har\-monic on $\Omega_c$, hence $\Cal
F$-continuous there. The set
 $$
\{z\in\Omega: d<f(z)<c\}=\{z\in\Omega_c: d<\max\{f(z),d\}<c\}
 $$
is therefore $\Cal F$-open, and hence $f$ is
$\Cal F$-continuous.

For (d) we refer to the proof given in \cite{EW3, Theorem 4.1}.
 {\hfill$\square$}\enddemo

\proclaim{Proposition 2.5} Every strongly $\Cal F$-pluri\-sub\-har\-monic
function $f:\Omega\to[-\infty,+\infty[\,$ is weakly $\Cal
F$-pluri\-sub\-har\-monic.
 \endproclaim

\demo{Proof} We may assume that $f$ is even $\Cal F$-cpsh. Let $(f_
\nu)$ be a sequence of finite con\-tinu\-ous pluri\-sub\-har\-monic
functions on open sets $\Omega_\nu$ containing $K$ from Definition
2.2(i) such that $f_\nu|K\to f|K$ uniformly. For any complex line $L$
in $\Bbb C^n$, $f_\nu|L\cap K'$ is finely hypoharmonic.  This uses
\cite{F1, Theorem 8.7} and the fact that the intersection of any $\Cal
F$-open subset of $\Bbb C^n$ with any complex line $L$ is finely open,
by Lemma 2.1. It follows by \cite{F1, Lemma 9.6} that $f|L\cap K'$ is
finely hypoharmonic, and in particular finely con\-tinu\-ous, by
\cite{F1, Theorem 9.10}. Consequently, $f$ is indeed weakly $\Cal
F$-pluri\-sub\-har\-monic.
 {\hfill$\square$}\enddemo

We now pass to concepts of $\Cal F$-holomorphic functions.  For a
compact set $K\subset\Bbb C^n$ we denote by $H_0(K)$ the algebra of
all restrictions to $K$ of ho\-lo\-morph\-ic functions defined on open
subsets of $\Bbb C^n$ containing $K$, and by $H(K)$ the closure of
$H_0(K)$ in $C(K,\Bbb C)$ (the continuous functions $K\to\Bbb C$ with
the uniform norm); then $H(K)$ is likewise an algebra.

\definition{Definition 2.6} Let $\Omega$ denote an $\Cal F$-open subset
of $\Bbb C^n$.

(i) (cf.\ \cite{EW2, Definition 6.1}.) A function $f:\Omega\to\Bbb C$
is said to be {\it strongly $\Cal F$-ho\-lo\-morph\-ic} if every point
of $\Omega$ has a compact $\Cal F$-neighborhood $K$ in $\Omega$ such
that $f|K\in{H(K)}$.

(ii) A function $f:\Omega\to\Bbb C$ is said to be {\it weakly $\Cal
F$-ho\-lo\-morph\-ic} if $f$ is $\Cal F$-con\-tinu\-ous and if, for
every complex line $L$ in $\Bbb C^n$, the restriction $f|L\cap\Omega$
is {\it finely ho\-lo\-morph\-ic}.
 \enddefinition

For finely ho\-lo\-morph\-ic functions see \cite{F3}, \cite{F7}.  Either
of the concepts strongly and weakly $\Cal F$-ho\-lo\-morph\-ic is an
$\Cal F$-local one. The class of all strongly, resp.\ weakly, $\Cal
F$-ho\-lo\-morph\-ic functions on $\Omega$ is an algebra, and the
latter class is closed under $\Cal F$-locally uniform convergence, in
view of \cite{F3, Th\'eor\graveaccent eme 4}. Clearly, every strongly
$\Cal F$-holo\-morph\-ic function is $\Cal F$-continuous (on $K'$ from
Definition 2.6(i), and so on all of $\Omega$).

If $f$ is strongly, resp.\ weakly, $\Cal F$-holomorphic on $\Omega$
($\Cal F$-open in $\Bbb C^n$) then the restriction of $f$ to
$L\cap\Omega$ ($L$ a $\Bbb C$-linear subspace $L\cong\Bbb C^k$ of
$\Bbb C^n$) has the same property on $L\cap\Omega$. This follows
easily from Lemma 2.1 above.

For $n=1$, $f$ is strongly (resp.\ weakly) $\Cal F$-ho\-lo\-morph\-ic
on $\Omega$ (finely open in $\Bbb C$) if and only if $f$ is finely
ho\-lo\-morph\-ic on $\Omega$. This is obvious in the weak case. In
the strong case, suppose first that $f$ is finely ho\-lo\-morph\-ic on
$\Omega$. By \cite{F3, Corollary, p.\,75}, every point of $\Omega$ has
a compact fine neighborhood $K$ in $\Omega$ such that $f|K\in R(K)$
($=H(K)$ in the $1$-dimensional case). Consequently, $f$ is indeed
strongly $\Cal F$-ho\-lo\-morph\-ic on $\Omega$. Conversely, if $f$ is
strongly $\Cal F$-holo\-morph\-ic then, for any compact set
$K\subset\Bbb C$, every function of class $H(K)$ is finely
holo\-morph\-ic on $K'$, see \cite{F3, p.\,63}.  With $K$ as in
Definition 2.6(i) this shows that $f$ indeed is finely holo\-morph\-ic on
$\Omega$.

\proclaim{Proposition 2.7} Every strongly $\Cal F$-ho\-lo\-morph\-ic
function $f:\Omega\to\Bbb C$ is weakly $\Cal F$-ho\-lo\-morph\-ic, and
in particular $\Cal F$-con\-tinu\-ous.
 \endproclaim

\demo{Proof} For any $K$ as in Definition 2.6(i) there exists a
sequence of holo\-morph\-ic functions $f_\nu$ defined on open sets
containing $K$ such that $f_\nu|K\to f|K$ uniformly.  For every
complex line $L$ in $\Bbb C^n$ this shows that the finely
ho\-lo\-morph\-ic functions $f_\nu|L\cap K'$ converge uniformly to
$f|L\cap K'$, which therefore is finely ho\-lo\-morph\-ic, see again
\cite{F3, p.\,63}. Consequently $f|L\cap\Omega$ is finely
ho\-lo\-morph\-ic, and so $f$ is indeed weakly $\Cal
F$-ho\-lo\-morph\-ic, being also $\Cal F$-con\-tinu\-ous.
{\hfill$\square$}\enddemo

The concept of weakly $\Cal F$-holo\-morph\-ic function can be
characterized in terms of weakly {\it $\Cal F$-pluri\-har\-monic}
functions (that is, functions $f:\Omega\to\Bbb C$ such that $\pm\Re f$
and $\pm\Im f$ are weakly $\Cal F$-pluri\-sub\-har\-monic on the $\Cal
F$-open set $\Omega\subset\Bbb C^n$):

\proclaim{Lemma 2.8} A function $f:\Omega\to\Bbb C$ is weakly $\Cal
F$-holo\-morph\-ic if and only if $f$ and each of the functions $z\mapsto
z_jf(z)$ $(j\in\{1,\ldots,n\})$ are weakly $\Cal
F$-pluri\-sub\-har\-monic on $\Omega$.
 \endproclaim

\demo{Proof} This reduces right away to the case $n=1$ which is due to
Lyons \cite{Ly}, cf.\ \cite{F3, Section 3}, and which asserts that a
function $h:U\to\Bbb C$, defined on a finely open set $U\subset\Bbb C$,
is finely holomorphic if and only if $h$ and $z\mapsto zh(z)$ are
(complex) finely har\-monic. 
 {\hfill$\square$}\enddemo

For any $\Cal F$-open set $U\subset\Bbb C^m$, an $n$-tuple
$(h_1,\dots,h_n)$ of strongly (resp.\ weakly) $\Cal F$-holomorphic
functions $h_j:U\to\Bbb C$ will be termed a strongly (resp.\ weakly)
$\Cal F$-{\it holomorphic map} $U\to\Bbb C^n$.
  
Assertion (b) of the following proposition provides two slight
strengthenings of \cite{EW2, Lemma 6.2}.

\proclaim{Proposition 2.9} Let $U\subset\Bbb C^m$ be $\Cal F$-open and
let $h=(h_1,\dots,h_n):U\to\Bbb C^n$ be a strongly $($resp.\ weakly$)$
$\Cal F$-holomorphic map.

{\rm (a)} The map $h:U\to\Bbb C^n$ is con\-tinu\-ous from $U$ with the
$\Cal F$-topology on $\Bbb C^m$ to $\Bbb C^n$ with the Euclidean
topology.

{\rm (b)} For any pluri\-sub\-har\-monic function $f$ on an open set
$\Omega$ in $\Bbb C^n$, the function $f\circ h$ is strongly $($resp.\
weakly$)$ $\Cal F$-pluri\-sub\-har\-monic on the $\Cal F$-open set
$h^{-1}(\Omega)=\{z\in U:h(z)\in\Omega\}\subset\Bbb C^m$.

{\rm (c)} For any ho\-lo\-morph\-ic function $f$ on an open set
$\Omega$ in $\Bbb C^n$, the function $f\circ h$ is strongly $($resp.\
weakly$)$ $\Cal F$-ho\-lo\-morph\-ic on the $\Cal F$-open set
$h^{-1}(\Omega)\subset\Bbb C^m$.
 \endproclaim

\demo{Proof} Assertion (a) holds because each $h_j$ (whether strongly
or weakly $\Cal F$-holo\-morph\-ic) is $\Cal F$-con\-tinu\-ous and
that the Euclidean topology on $\Bbb C^n$ is the product of the
Euclidean topology on each of $n$ copies of $\Bbb C$.

For (b) with each $h_j$ {\it strongly} $\Cal F$-holomorphic we begin
by showing that, if the pluri\-sub\-har\-monic function $f$ on
$\Omega$ is finite and continuous, then $f\circ h$ is even $\Cal
F$-{\it cpsh} (cf.\ Definition 2.2(i)) on $h^{-1}(\Omega)$, which is
$\Cal F$-open according to (a).  Every point $a\in h^{-1}(\Omega)$ has
a compact $\Cal F$-neighborhood $K_j$ in $h^{-1}(\Omega)$ ($\subset
U\subset\Bbb C^m$) such that $h_j|K_j\in H(K_j)$. Thus there exists a
sequence $(h_j^\nu)_{\nu\in\Bbb N}$ of holo\-morph\-ic functions
$h_j^\nu$ on open sets $U_j^\nu$ in $\Bbb C^m$ containing $K_j$ such
that $h_j^\nu|K_j\to h_j|K_j$ uniformly as $\nu\to\infty$.  Write
$K=K_1\cap\ldots\cap K_n$ and $h^\nu=(h_1^\nu,\dots,h_n^\nu)$ on
$U^\nu=U_1^\nu\cap\ldots\cap U_n^\nu$. Then $h_j|K_j\in C(K_j,\Bbb C)$
and hence $h|K\in C(K,\Bbb C^n)$. It follows that $h(K)$ is a compact
subset of $\Omega\subset\Bbb C^n$.  Denoting by $\|\,\cdot\,\|$ the
Euclidean norm on $\Bbb C^n$ and by $B$ the closed unit ball in $\Bbb
C^n$, there exists accordingly $\delta>0$ such that $h(K)+\delta
B\subset\Omega$. We may assume that
$\|h^\nu(z)-h(z)\|\allowmathbreak<\delta$ for any $\nu$ and any $z\in
K$.  Under the present extra hypothesis, $f$ is finite and uniformly
continuous on the compact set $h(K)+\delta B$ containing any
$h^\nu(K)$, and it follows that $f\circ h^\nu|K\to f\circ h|K$
uniformly as $\nu\to\infty$. Because $f\circ h^\nu$ is finite,
con\-tinu\-ous, and pluri\-sub\-har\-monic, on the open set
$U^\nu\supset K$, we have $f\circ h|K\in S(K)$. By varying $a\in
h^{-1}(\Omega)$ and hence the $\Cal F$-neighborhood $K$ of $a$ in
$h^{-1}(\Omega)$ we infer that $f\circ h$ is $\Cal F$-cpsh on
$h^{-1}(\Omega)$.

If we drop the extra hypothesis that $f$ be finite and continuous, $f$
is the pointwise limit of a decreasing net of finite continuous
pluri\-sub\-har\-monic functions $f_\nu$ on $\Omega$, and $f\circ h$
is then the pointwise limit of the decreasing net of functions
$f_\nu\circ h$ on $h^{-1}(\Omega)$ which we have just shown are
$\Cal F$-cpsh, and so $f\circ h$ is indeed strongly $\Cal
F$-pluri\-sub\-har\-monic, cf.\ Definition 2.2(ii) (with $\Omega$
replaced by $h^{-1}(\Omega)$).

Next suppose instead that each $h_j$ is {\it weakly} $\Cal
F$-holo\-morph\-ic on $U$, and consider a complex line $L$ in $\Bbb
C^m$; then $L\cap U$ is finely open in $L$.  According to Definition
2.6(ii), $h_j|L\cap U$ is then finely holo\-morph\-ic, which is the
same as strongly $\Cal F$-holo\-morph\-ic (see above for $n=1$). As
shown above (now with $m=1$ and with $U$ replaced by $L\cap U$) it
follows that $f\circ h|L\cap h^{-1}(\Omega)$ is strongly $\Cal
F$-pluri\-sub\-har\-monic, which is the same as finely hypoharmonic
(because the dimension is $1$). According to Definition 2.2(iii) this
means that $f\circ h$ indeed is weakly $\Cal F$-pluri\-sub\-har\-monic
on $h^{-1}(\Omega)$, noting that $f\circ h$ is $\Cal F$-upper
semi\-con\-tinu\-ous in view of (a) because $f$ is upper
semi\-con\-tinu\-ous.

For (c), suppose first that each $h_j$ is {\it strongly} $\Cal
F$-holo\-morph\-ic on $U$.  Proceeding as in the first part of the
proof of (b) we arrange that $f\circ h^\nu|K\to f\circ h|K$ uniformly
as $\nu\to\infty$; but now $f\circ h^\nu$ is holo\-morph\-ic on
$U^\nu$. We therefore conclude that $f\circ h|K\in{H(K)}$, and so
$f\circ h$ is indeed strongly $\Cal F$-holo\-morph\-ic according to
Definition 2.6(i).

If instead each $h_j$ is {\it weakly} $\Cal F$-holo\-morph\-ic on $U$
then, for every complex line $L$ in $\Bbb C^m$, each $h_j|L\cap U$ is
again {\it strongly} $\Cal F$-holo\-morph\-ic. As just established,
this implies that $f\circ h|L\cap h^{-1}(\Omega)$ is strongly $\Cal
F$-holo\-morph\-ic, or equivalently finely holo\-morph\-ic. We
conclude that indeed $f\circ h$ is {\it weakly} $\Cal
F$-holo\-morph\-ic, according to Definition 2.6(ii), noting that $f\circ
h$ is $\Cal F$-con\-tinu\-ous in view of (a).
{\hfill$\square$}\enddemo

In the version of Proposition 2.9 with \lq weakly' in each of the
three occurrences one may allow $f$ in (b) to be just {\it weakly
$\Cal F$-pluri\-sub\-har\-monic} (in place of pluri\-sub\-har\-monic),
and similarly $f$ in (c) to be {\it weakly $\Cal F$-holo\-morph\-ic}
(in place of holo\-morph\-ic), see Theorem 4.6 at the end of the
paper.  At this point we merely show that we may allow $f$ in (b) (of
Proposition 2.9) to be {\it strongly} $\Cal F$-pluri\-sub\-har\-monic,
and $f$ in (c) to be strongly $\Cal F$-holo\-morph\-ic:

\proclaim{Theorem 2.10} Let $U\subset\Bbb C^m$ be $\Cal F$-open and let
$h=(h_1,\dots,h_n):U\to\Bbb C^n$ be a weakly $\Cal F$-holomorphic map.

{\rm (a)} The map $h:U\to\Bbb C^n$ is con\-tinu\-ous from $U$ with the
$\Cal F$-topology on $\Bbb C^m$ to $\Bbb C^n$ with the $\Cal
F$-topology there.

{\rm (b)} For any strongly $\Cal F$-pluri\-sub\-har\-monic function
$f$ defined on an $\Cal F$-open set $\Omega$ in $\Bbb C^n$, the
function $f\circ h$ is weakly $\Cal F$-pluri\-sub\-har\-monic on the
$\Cal F$-open set $h^{-1}(\Omega)=\{z\in U:h(z)\in\Omega\}\subset\Bbb
C^m$.

{\rm (c)} For any strongly $\Cal F$-ho\-lo\-morph\-ic function $f$
defined on an $\Cal F$-open set $\Omega$ in $\Bbb C^n$, the function
$f\circ h$ is weakly $\Cal F$-ho\-lo\-morph\-ic on the $\Cal F$-open
set $h^{-1}(\Omega)\subset\Bbb C^m$.
 \endproclaim
 
\demo{Proof} For the present weakly $\Cal F$-holo\-morph\-ic functions
$h_j$ assertion (a) is stronger than Proposition 2.9(a).  We shall
prove that $h^{-1}(\Omega)$ is $\Cal F$-open in $\Bbb C^m$ for any
$\Cal F$-open set $\Omega$ in $\Bbb C^n$.  Fix a point $a\in
h^{-1}(\Omega)$ and write $h(a)=b$ ($\in\Omega$). According to
\cite{EW2, Theorem 2.3} there exist a pluri\-sub\-har\-monic function
$\phi$ on an open ball $B(b,r)$ in $\Bbb C^n$ and a number $c<\phi(b)$
such that the basic $\Cal F$-neighborhood
 $$ 
W=\{w\in B(b,r):\phi(w)>c\}
 $$ 
of $b$ in $\Bbb C^n$ is a subset of the $\Cal F$-open set $\Omega$ in
$\Bbb C^n$. Then $h^{-1}(W)\subset h^{-1}(\Omega)$ ($\subset U$), and
 $$ 
h^{-1}(W)=\{z\in U:h(z)\in B(b,r)\text{ \,and \,}(\phi\circ
 h)(z)>c\}
 $$
is $\Cal F$-open in $\Bbb C^m$ because $h:U\to\Bbb C^n$ is $\Cal
F$-continuous by Proposition 2.9(a) and that $\phi\circ
h:h^{-1}(B(b,r))\to[-\infty,+\infty[\,$ is weakly $\Cal
F$-pluri\-sub\-har\-monic by Proposition 2.9(b) (applied with
$\Omega,f$ replaced by $B(b,r),\phi$), and in particular $\Cal F$-{\it
continuous}, by Theorem 2.4(c).  By varying $a\in h^{-1}(\Omega)$ we
infer that indeed $h^{-1}(\Omega)$ is $\Cal F$-open.

For (b) we may assume that $f$ is even $\Cal F$-{\it cpsh} on
$\Omega$. Let $K\subset\Omega$ be as in Definition 2.2(i), and let
$(f^\nu)$ be a sequence of finite con\-tinu\-ous
pluri\-sub\-har\-monic functions on open sets $\Omega^\nu\supset K$
such that $f^\nu|K\to f|K$ uniformly as $\nu\to\infty$. According to
Proposition 2.9(a),(b) each $f^\nu\circ h$ is weakly $\Cal
F$-pluri\-sub\-har\-monic on the $\Cal F$-open set
$h^{-1}(\Omega^\nu)$.  By~(a), $h^{-1}(K')$ is $\Cal F$-open, and it
follows that each $f^\nu\circ h|h^{-1}(K')$ likewise is weakly $\Cal
F$-pluri\-sub\-har\-monic, in particular $\Cal F$-upper
semicontinuous; and hence so is its uniform limit $f\circ
h|h^{-1}(K')$ in view of \cite{F1, Lemma 9.6}.  By varying
$K\subset\Omega$ and hence $K'$ we conclude that indeed $f\circ h$ is
weakly $\Cal F$-pluri\-sub\-har\-monic on $h^{-1}(\Omega)$.

Finally, the proof of (c) is quite parallel to that of (b) in view of
Proposition 2.9(a),(c), using \cite{F3, Th\'eor\graveaccent eme 4} in
place of \cite{F1, Lemma 9.6}. 
 {\hfill$\square$}\enddemo

Theorem 2.10 has two corollaries for $m=1$ and $n=1$, respectively. In
either corollary \lq strongly' can be replaced by \lq weakly'
according to Theorem 4.6. For $m=1$ we have

\proclaim{Corollary 2.11} Let $h_j:U\to\Bbb C$ $(j\in\{1,\dots,n\})$
be finely holomophic functions defined on a finely open set
$U\subset\Bbb C$, and write $h=(h_1,\dots,h_n)$. For any strongly
$\Cal F$-pluri\-sub\-har\-monic $($resp.\ strongly $\Cal
F$-holo\-morph\-ic$)$ function $f$ defined on an $\Cal F$-open set
$\Omega\subset\Bbb C^n$, the function $f\circ h$ is finely
hypo\-harm\-onic $($resp.\ finely holo\-morph\-ic$)$ on the finely
open set $h^{-1}(\Omega)=\{z\in U:h(z)\in\Omega\}\subset\Bbb C$.
 \endproclaim

\demo{Remark 2.12} If it can be proved that the cone of $\Cal F$-cpsh,
resp. the algebra of strongly $\Cal F$-holo\-morph\-ic, functions on
an $\Cal F$-open subset $\Omega$ of $\Bbb C^n$ is closed under uniform
convergence, then the proofs of Theorem 2.10(b),(c) easily show that
\lq weakly' can be replaced throughout the theorem by \lq
strongly'. Indeed, with $f$ $\Cal F$-cpsh (resp.\ strongly $\Cal
F$-holo\-morph\-ic) and $h$ strongly $\Cal F$-holo\-morph\-ic, let
$K\subset\Omega$ denote a compact $\Cal F$-neighborhood of a point
$a\in\Omega$, and let $(f_\nu)$ denote a sequence of finite continuous
pluri\-sub\-har\-monic (resp.\ a sequence of holo\-morph\-ic)
functions, defined on open subsets $\Omega_\nu$ of $\Bbb C^n$
containing $K$, and such that $f_\nu\to f$ uniformly on $K$. Then
$f_\nu\circ h$ is strongly $\Cal F$-pluri\-sub\-har\-monic and even
$\Cal F$-cpsh (resp.\ strongly $\Cal F$-holo\-morph\-ic) on the $\Cal
F$-open set $h^{-1}(K')\subset h^{-1}(\Omega_\nu)\subset U\subset\Bbb
C^m$, by Proposition 2.9 (cf.\ the beginning of the proof of Theorem
2.9(b)).  Under the stated extra hypothesis it will follow that the
uniform limit $f\circ h$ of $(f_\nu\circ h)$ on $h^{-1}(K')$ likewise
is $\Cal F$-cpsh (resp.\ strongly $\Cal F$-holo\-morph\-ic) on
$h^{-1}(K')$, and therefore on $h^{-1}(\Omega)$, by varying $K$ and
hence $K'$.---If $f$ is merely strongly $\Cal F$-{\it
pluri\-sub\-har\-monic} (rather than $\Cal F$-cpsh), it follows as
usual that indeed $f\circ h$ is likewise strongly $\Cal
F$-pluri\-sub\-har\-monic.
 \enddemo

In the case $n=1$ the extra hypothesis stated in the above remark is
always fulfilled in view of \cite{F1, Lemma 9.6} (resp.\ \cite{F3,
Th\'eor\graveaccent eme 4}). We therefore have the following corollary
of Theorem 2.10 for that case:

\proclaim{Corollary 2.13} Let $h:U\to\Bbb C$ be a strongly $\Cal
F$-holo\-morph\-ic function defined on an $\Cal F$-open set
$U\subset\Bbb C^m$. For any finely hypo\-har\-mon\-ic $($resp.\ finely
holo\-morph\-ic$)$ function $f$ defined on a finely open set
$\Omega\subset\Bbb C$ the function $f\circ h$ is strongly $\Cal
F$-pluri\-sub\-har\-monic $($resp.\ strongly $\Cal F$-holo\-morph\-ic$)$
on the $\Cal F$-open set $h^{-1}(\Omega)\subset\Bbb C^m$.
 \endproclaim

The same with \lq strongly' replaced throughout by \lq weakly' is
simply the case $n=1$ of Theorem 2.10 as it stands.

\proclaim{Proposition 2.14} Let now $\Omega$ be a Euclidean open
subset of \,$\Bbb C^n$. For a function
$f:\Omega\to[-\infty,+\infty[\,$ the following are equivalent:

\roster
\item"{\rm(i)}" $f$ is pluri\-sub\-har\-monic $($in the ordinary sense$)$.

\item"{\rm(ii)}" $f$ is strongly $\Cal F$-pluri\-sub\-har\-monic and not
identically $-\infty$ on any component of \,$\Omega$. 

\item"{\rm(iii)}" $f$ is weakly $\Cal F$-pluri\-sub\-har\-monic and not
identically $-\infty$ on any com\-po\-nent of \,$\Omega$.
 \endroster\endproclaim
 
\demo{Proof} Every finite continuous pluri\-sub\-har\-monic function
on $\Omega$ is of course $\Cal F$-cpsh.  It follows that any
pluri\-sub\-har\-monic function on $\Omega$ is strongly $\Cal
F$-pluri\-sub\-har\-monic (being the pointwise limit of a decreasing
sequence of finite continuous pluri\-sub\-har\-monic
functions). Conversely, if $f$ is weakly $\Cal
F$-pluri\-sub\-har\-monic on $\Omega$ then $f$ is
pluri\-sub\-har\-monic on every connectivity component $\omega$ of
$\Omega$ on which $f$ is not identically $-\infty$. To see this, first
observe that $f$ is $\Cal F$-locally bounded from above, so every
point $a\in\Omega$ has an $\Cal F$-neighborhood $U\subset\Omega$ on
which $f<M_a$, say.  According to \cite{EW2, Proposition 4.1} we may
further arrange that there exists $\delta>0$ such that, for every
complex line $L$ passing through $a$, the intersection $U\cap L$
contains a circle about $a$ with radius at least $\delta$. By the
maximum principle for finely subharmonic functions on a planar domain,
\cite{F2, Theorem 2.3}, it follows that $f<M_a$ on the discs bounded
by these circles, hence in particular on the ball
$B(a,\delta_a)$.---For functions that are locally bounded from above
in the Euclidean topology, the statement alternatively follows from
\cite{F1, Theorem 9.8(a)} in view of \cite{Le1, D\'efinition 1 (p.\
306)}.  {\hfill$\square$}\enddemo

\demo{Remark 2.15} Similarly to Proposition 2.14, a function
$f:\Omega\to\Bbb C$ (with $\Omega$ Euclidean open) is holo\-morph\-ic
if and only if $f$ is strongly, or equivalently weakly, $\Cal
F$-holomorphic; the `if part' follows from \cite{F3, p.\,63} in view
of Hartogs' theorem.
 \enddemo

We close this section with an application to pluripolar hulls.  Recall
that the {\it pluripolar hull} $P^*_\Omega$ of a pluripolar set $P$
relative to an open set $\Omega$ containing $P$ is defined as the
following set (closed relatively to $\Omega$):
 $$
P^*_\Omega=\bigcap_u\{z\in\Omega:u(z)=-\infty\},
 $$
where the intersection is taken over all pluri\-sub\-har\-monic functions
$u$ defined on $\Omega$ and such that $u|P\equiv-\infty$. A
pluri\-polar set $E$ has empty plurifine interior $E'$, \cite{EW2,
Theorem 5.2}. (More generally, a polar set is a Lebesgue null set and
has therefore empty fine interior.) 
 
For any set $E\subset\Bbb C^m$, $m\in\Bbb N$, and any function
$h:E\to\Bbb C$ we denote by $\Gamma_h(E)=\{(z,h(z)):z\in E\}$ the
graph of $h|E$ and by $\Gamma_h(E)^*_{\Bbb C^{m+1}}$ the pluripolar
hull of $\Gamma_h(E)$.

\proclaim{Proposition 2.16} Let $h$ be a weakly $\Cal F$-holomorphic
function on an $\Cal F$-domain $U\subset\Bbb C^m$.

{\rm(a)} If \,$h\not\equiv0$, the set $h^{-1}(0)$ of zeros of \,$h$ is
pluripolar in $\Bbb C^m$. In particular, the graph $\Gamma_h(U)$ of
$h$ is pluripolar in $\Bbb C^{m+1}$.

{\rm(b)} If \,$E$ is a non-pluripolar subset of \,$U$ then $\Gamma_h(U)$
is pluripolar, and $\Gamma_h(U)\subset\Gamma_h(E)^*_{\Bbb
C^{m+1}}$.
 \endproclaim

With $h$ supposed {\it strongly} $\Cal F$-holomorphic on $U$,
Proposition was obtained in \cite{EW3, Corollary 4.4 and Theorem 4.5},
extending \cite{EW2, Theorem 6.4}, and \cite{EEW, Theorem 3.5}

\demo{Proof of Proposition 2.16} (a) It follows from Proposition
2.9(b) that the function $\log|h(\,\cdot\,)|$ is weakly $\Cal
F$-pluri\-sub\-har\-monic on $U$. Since $\log|h(z)|=-\infty$ for $z\in
h^{-1}(0)$, but $\log|h(\,\cdot\,)|\not\equiv-\infty$, we conclude
from Theorem 2.4(d) that $h^{-1}(0)$ is pluripolar.

(b) The function $(z,w)\mapsto w-h(z)$ is weakly $\Cal F$-holomorphic
and $\not\equiv0$ on the $\Cal F$-open set $U\times\Bbb C\subset\Bbb
C^{m+1}$. Again by Proposition 2.9(b) it follows that the function
$(z,w)\mapsto\log| w-h(z)|$ is weakly $\Cal F$-pluri\-sub\-har\-monic
and $\not\equiv-\infty$ on $U\times\Bbb C$. Since this function equals
$-\infty$ on $\Gamma_h(E)$ we conclude that $\Gamma_h(E)$ is
pluripolar. By Josefson's theorem \cite{J} there exists a
pluri\-sub\-har\-monic function $f$ on all of $\Bbb C^{m+1}$ such that
$f(z,h(z))=-\infty$ for every $z\in V$. It follows by Theorem 2.3(a)
that $f(z,h(z))=-\infty$ even for every $z\in U$, and hence
$\Gamma_h(U)$ is pluripolar in $\Bbb C^{m+1}$. By definition of the
pluripolar hull of $\Gamma_h(E)$ we conclude that indeed
$\Gamma_h(U)\subset\Gamma_h(E)^*_{\Bbb C^{m+1}}$.
{\hfill$\square$}\enddemo

\heading{\bf 3. A characterization of weakly $\Cal
F$-pluri\-sub\-har\-monic functions}\endheading

By the prefix `$\Bbb R^{2n}$-fine' we denote concepts relative to the
Cartan fine topology on $\Bbb C^n\cong\Bbb R^{2n}$. Recall that this
topology is finer than the plurifine topology $\Cal F$, \cite{F6}.  It
is well-known that a plurisubharmonic function $f$ on a domain
$\Omega\subset\Bbb C^n$ is subharmonic when considered as a function
on $\Omega\subset\Bbb R^{2n}$, because the average of $f$ over a
sphere can be expressed in terms of the average of $f$ over the
circles that are intersection of the sphere with complex lines passing
through the center. While this approach does not work in the fine
setting, the analogous result nevertheless remains valid.  Indeed, a
well-known characterization of pluri\-sub\-har\-monic functions (see
\cite{Le2, Th\'eor\graveaccent eme 1 (p.\ 18)} or \cite{K, Theorem
2.9.12}) may be adapted as follows. This will lead to further
properties of weakly $\Cal F$-pluri\-sub\-har\-monic functions
(Theorems 3.7 and 3.9).

\proclaim{Theorem 3.1} A function $f:\Omega\to[-\infty,+\infty[$
$(\Omega$ $\Cal F$-open in $\Bbb C^n)$ is weakly $\Cal
F$-pluri\-sub\-har\-monic if and only if $f\circ h$ is $\Bbb
R^{2n}$-finely hypoharmonic on the $\Cal F$-open set $h^{-1}(\Omega)$
for every $\Bbb C$-affine bijection $h$ of \,$\Bbb C^n$.
 \endproclaim

For the proof of the `only if part' of Theorem 3.1 we need the following

\proclaim {Lemma 3.2} Let $u_1,u_2$ be bounded subharmonic functions
on an open set $B\subset\Bbb
R^n$, and consider the function $f=u_1-u_2$ on $B$. Let $U$ be a
finely open Borel subset of $B$. Then $f|U$ is finely subharmonic if and
only if the signed Riesz measure $\Delta f$ on $B$ has a positive
restriction to $U$.
 \endproclaim
 
\demo{Proof} Suppose first that $(\Delta f)|U\ge0$, and let us prove
that $f$ then is finely subharmonic on $U$.

Recall that the \rm{base} $b(X)$ in $B$ of $X\subset B$ consists of
the points of $B$ at which $X$ is not thin. Denote by $\complement X$
the complement of $X$ relative to $B$. For a finely open set $U$, its
regularization equals $r(U)=\complement b(\complement U))= U\cup
i(\complement U)$ where $i(X)$ is the polar set consisting of the
points of $X $ at which $X$ is thin. We may assume that $U$ is {\it
regular}, i.e., $U=r(U)$ and hence an $F_\sigma$-set, for $u_1$ and
$u_2$ are bounded, and $\Delta u_1$ and $\Delta u_2$ therefore do not
charge the polar set by which $U$ differs from $r(U)$.

Writing $\Delta f:=\mu=\mu_+-\mu_-$ on $B$ we have, by hypothesis,
$\mu_-|U=0$. Proceeding as in the proof of the former (and easier) `if
part' of \cite{F1, Theorem 8.10}, consider any bounded finely open set
$V$ of compact closure $\overline V\subset U$. Then
$(\mu_-)^{\complement V}=\mu_-$ because $\mu_-$ is carried by
$b(\complement V)\supset\complement U$. For any $x\in V$ we obtain in
terms of the Green kernel $G$ on $B$ according to \cite{D, Theorems
1.X.3 and 1.X.5} applied within $B$
 $$\align 
\int G\mu_-d\epsilon_x^{\complement V}&=\wh R_{G\mu_-}^{\complement
V}(x) =G\bigl((\mu_-)^{\complement V}\bigr)(x)=G\mu_-(x)<+\infty,\\
\int G\mu_+d\epsilon_x^{\complement V}&=\wh R_{G\mu_+}^{\complement
V}(x) \le G\mu_+(x)<+\infty,
 \endalign$$
whence by subtraction $\int G\mu\,d\epsilon _x^{\complement V} \le
G\mu(x)$, showing that the finely continuous function $G\mu$ is finely
hyperharmonic, and indeed (being also bounded) finely superharmonic on
$U$, \cite{F1, Theorem 8.10 and \S10.4}. By the Riesz representation
theorem, $f=-G\mu+h$ on $B$, with $h$ harmonic on $B$, in particular
finely harmonic on $U$, and $f$ is therefore likewise finely
subharmonic on $U$.

Conversely, suppose that $f|U$ is finely subharmonic.  Recall the
corollary in \cite{D, 1.XI.18} that if two subharmonic
functions $g_1$ and $g_2$, defined on some open set, coincide on a set
$A$, then their Riesz masses satisfy $\Delta g_1=\Delta g_2$ on the
fine interior of $A$.  Hence the Riesz measure $\Delta f=\Delta
u_1-\Delta u_2$ on $U$ is independent of the choice of $u_1$ and
$u_2$. In fact, if $f=w_1-w_2$ (with $w_1$ and $w_2$ subharmonic on
$U$) then $u_1+w_2=u_2+w_1$ on $U$, hence $\Delta u_1+\Delta
w_2=\Delta u_2+\Delta w_1$ on (the fine interior of) $U$, that is,
$\Delta w_1-\Delta w_2=\Delta u_1-\Delta u_2$ on $U$.

Since $f$ is finely subharmonic on $U$ it follows by the proof of
\cite{F1, Theorem 9.9} that every point $x\in U$ has a fine
neighborhood $V_x\Subset U$ in which we can write $f=v_1-v_2$, where
$v_1$ and $v_2$ are superharmonic functions on some Euclidean
neighborhood $B_0$ of $x$ in $B$.  Moreover, $v_2$ is the swept-out on
$B_0\setminus V_x$ of a certain superharmonic function $\ge0$ on
$B_0$.  The Riesz mass of $v_2$ is concentrated on the fine boundary
of the complement of $V_x$, cf.\ e.g.\ \cite{D, Theorem 1.XI.14(b)},
hence on the fine boundary of $V_x$.  It follows that the Riesz mass
$\Delta f$ of $f$ is positive on $V_x$ for every $x$. By the
quasi-Lindel\"of property, we can find countably many $x_j\in U$ such
that $U=\bigcup_{j=1}^\infty V_{x_j}\cup E$, where $E$ is
polar. Clearly $\Delta f$ is positive on $\bigcup
V_{x_i}$. Because $E$ is polar, the Riesz mass of a bounded
subharmonic function does not charge $E$, so we have $\Delta u_i(E)=0$
($i=1,2$). We conclude that the measure $\Delta f$ is positive on $U$.
{\hfill$\square$}\enddemo

\demo{Proof of the \lq only if part' of Theorem 3.1} The proof merely
uses that $f$ is finely hypoharmonic in each variable separately.
According to Theorem 2.4(a), every point $a\in\Omega$ has a bounded
$\Cal F$-open $\Cal F$-neighborhood $O\subset\Omega$ on which $f$ is
representable as $f=u_1-u_2$, where $u_1$ and $u_2$ are bounded
plurisubharmonic functions, defined on an open ball $B$ in $\Bbb C^n$
containing $O$. For every $j\in\{1,\dots,n\}$ the distributions
 $$ 
M_{i,j}=M_{i,j}(z) =\frac{\partial^2u_i(z_1,\dots,z_n)}{\partial
z_j\partial\bar z_j},\qquad i=1,2,
 $$
are well defined positive measures on $B$. Below we show that
$(M_{1,j}-M_{2,j})|O\ge0$ and hence
$\Delta(u_1-u_2)=4\sum_{j=1}^n(M_{1,j}-M_{2,j})\ge0$ on $O$, where
$\Delta$ denotes the Laplacian on $\Bbb C^n\cong\Bbb R^{2n}$.
According to Lemma 3.2 this implies that $u_1-u_2$ indeed is $\Bbb
R^{2n}$-finely subharmonic on $O$, and hence actually on all of
$\Omega$, by varying $a$ and $O$.

For the proof that $(M_{1,j}-M_{2,j})|O\ge0$ it is convenient to write
points of $\Bbb C^n$ as $(z,w)$, now with $z\in\Bbb C$ and $w\in\Bbb
C^{n-1}$. Each of the above measures $M_{i,j}$ on $B$ then takes the
form
 $$
M_i=M_i(z,w) =\frac{\partial^2u_i(z,w)}{\partial z\partial\bar
z},\qquad i=1,2,
 $$ and we shall prove that $(M_1-M_2)|O\ge0$.

For each $w\in\Bbb C^{n-1}$ define
 $$ 
B(w)=\{z\in\Bbb C:(z,w)\in B\}
 $$ 
(and similarly with $B$ replaced by other subsets of $\Bbb C^n$). The
functions $u_i(\,\cdot\,,w)$, $i=1,2$, induce subharmonic functions
$u_i(\,\cdot\,,w)$ on the open subset $B(w)$ of $\Bbb C$, and the
distributions
 $$ 
\mu_{i,w}=\mu_{i,w}(z)= \frac{\partial^2u_i(z,w)}{\partial
z\partial\bar z},\qquad i=1,2,
 $$
are therefore positive measures on the open set $B(w)$ (if non-empty).
Being 
weakly $\Cal F$-plurisubharmonic on $O$,
$f=u_1-u_2$ induces the finely subharmonic function
$f(\,\cdot\,,w)=u_1(\,\cdot\,,w)-u_2(\,\cdot\,,w)$ on the finely open
set $O(w)$.  According to the planar version of Lemma 3.2, applied to
the induced bounded subharmonic functions $u_i(\,\cdot\,,w)$ on
$B(w)$, $i=1,2$, the Riesz measure $\mu_{1,w}-\mu_{2,w}$ of
$f(\,\cdot\,,w)$ is positive on the finely open set $O(w)\subset
B(w)$.

Let $V_z$, $V_w$ denote Lebesgue measure on $\Bbb C$, $\Bbb C^{n-1}$,
respectively. For any test function $\phi\in C_0^\infty(B)$ we have by
Fubini's theorem
 $$
\align \int_B\phi\,dM_i&=\int_B\frac{\partial^2\phi(z,w)}{\partial
z\partial\bar z} u_i(z,w)\,dV_z\,dV_w\\ &=\int_{\Bbb
C^{n-1}} \Bigl(\int_{B(w)}\frac{\partial^2\phi(z,w)}{\partial
z\partial\bar z} u_i(z,w)\,dV_z\Bigr)\,dV_w\\ &=\int_{\Bbb
C^{n-1}} \Bigl(\int_{B(w)}\phi(z,w)\,d\mu_{i,w}(z)\Bigr)\,dV_w.
 \endalign$$

Choose a compact $\Cal F$-neighborhood $K$ of the given point $a\in
O\subset\Omega$ so that $K\subset O$. There exists a decreasing
sequence of functions $\phi_k\in C_0^\infty(B)$ with $0\le\phi_k\le1$
so that $\phi_k=1$ on $K$ and $\phi_k\searrow\chi_K$ (the characteristic
function of $K$) as $k\nearrow\infty$. Since $M_i$ and $\mu_{i,w}$
are locally finite positive measures and $B$ and $\phi_k$ are bounded
we obtain by the monotone convergence theorem
 $$
\align M_i(K)&=\int_B\chi_K\,dM_i=\lim_{k\to\infty}\int_B\phi_k\,dM_i\\
&=\lim_{k\to\infty}\int_{\Bbb C^{n-1}}
\Bigl(\int_{B(w)}\phi_k(\,\cdot\,,w)\,d\mu_{i,w}\Bigr)\,dV_w\\
&=\int_{\Bbb C^{n-1}}
\Bigl(\int_{B(w)}\chi_{K(w)}\,d\mu_{i,w}\Bigr)\,dV_w\\
&=\int_{\Bbb C^{n-1}}\mu_{i,w}(K(w))\,dV_w.\\ \intertext{It
follows that} M_1(K)-M_2(K)&=\int_{\Bbb C^{n-1}}
\bigl(\mu_{1,w}(K(w))-\mu_{2,w}(K(w))\bigr)\,dV_w\ge0
 \endalign$$
because $\mu_{1,w}-\mu_{2,w}\ge0$ on $O(w)\supset K(w)$. Thus
$M_1(K)\ge M_2(K)$ for every compact $\Cal F$-neighborhood $K$ of $a$
in $O$.

The proof of Theorem 2.4(a) shows that we may take $O=\{z\in
B(z_0,r):\Phi^*(z)\ge-\frac14\}$, where $\Phi^*$ is plurisubharmonic
on the open ball $B(z_0,r)$, and in particular upper semicontinuous
there. It follows that $O=\bigcup_{p\in\Bbb N}F_p$ with
 $$ 
F_p=\{z\in B(z_0,(1-\tfrac1p)r):\Phi^*(z)\ge-\tfrac14+\tfrac1p\},
 $$ 
a bounded closed and hence compact subset of $\Bbb C^n$. Defining
$K_p=F_p\cup K$ we find that $K_p$ is a compact $\Cal F$-neighborhood
of $a$. We infer that $K_p\nearrow O$ as $p\nearrow\infty$, and
consequently
 $$ 
M_1(O)=\sup_{p\in\Bbb N} M_1(K_p)\ge\sup_{p\in\Bbb N} M_2(K_p)=M_2(O).
 $$ 
By Lemma 3.2, this completes the proof of the `only if part' of Theorem
3.1.  \hfill{$\square$}\enddemo

For the proof of the 'if part' of Theorem 3.1 we will need the
following lemma, and some results of Bedford and Taylor on slicing of
currents. 

\proclaim{Lemma 3.3} Let $f$ be a bounded finely subharmonic function
on an $\Cal F$-open set $\Omega\subset\CC^n$ and suppose that for
every $\Bbb C$-affine bijection $h$ of $\Bbb C^n$ the function $f\circ
h$ is finely subharmonic on $h^{-1}(\Omega)$. Then every $z_0\in
\Omega$ admits a $($compact$)$ $\Cal F$-neighborhood $K_{z_0}$ such that
$f$ can be written as
 $$
f=f_1-f_2\quad \text{on\ } K_{z_0},
 $$
where $f_1, f_2$ are plurisubharmonic functions defined on a ball
$B(z_0,r)\supset K_{z_0}$.
\endproclaim

\demo{Proof} As in the proof of (a) of Theorem 2.4. we can assume that
$-1<f<0$, and find a compact $\Cal F$-neighborhood $V$ of $z_0$ and a
negative pluri\-sub\-har\-monic function $\phi$ on a ball
$B(z_0,r)\supset V$ such that $\phi(z_0)=-1/2$ and $\phi= -1$ on
$B(z_0,r)\setminus V$.  For every $\lambda>0$ we can form the function
 $$
 u_\lambda(z)=\cases\max\{-\lambda, f(z)+\lambda\phi(z)\} &
\text{for\ } z\in \Omega\cap B(z_0,r), \\
 -\lambda& \text{for\ } z\in
B(z_0,r)\setminus V.
 \endcases $$
It is a bounded finely subharmonic function on $B(z_0,r)$, hence
$u_\lambda$ is subharmonic on $B(z_0,r)$.  Similarly, for every $\Bbb
C$-affine bijection $h$ of $\Bbb C^n$ the function $u_\lambda\circ h$
is finely subharmonic, hence subharmonic on $h^{-1}(B(z_0,r))$. From
this we conclude that $u_\lambda$ is in fact plurisubharmonic.  Taking
$\lambda=4$, we see that
 $$
u_4(z)=f(z)+4\phi(z)
 $$
on the closed $\Cal F$-neighborhood
$K_{z_0}=\{z\in\Omega:\phi(z)\ge-2/3\}\subset V\cap B(z_0,r)$, and
$K_{z_0}$ is compact along with $V$.  This proves the lemma.
\hfill{$\square$}\enddemo

\proclaim{Corollary 3.4} We keep the notation as above. Then for every
$z_0\in\Omega$, $f$ is $\Cal F$-continuous on $K_{z_0}$, hence on
$\Omega$.
\endproclaim

We recall from \cite{BT3} the concept of slice of an
$(n-1,n-1)$-current, now on a domain $D$ in $\Bbb C^n$. As usual we
will write $d=\partial +\overline{\partial}$ and
$d^c=i(\overline\partial -{\partial})$ so that $dd^c u= 2i
\partial\overline{\partial} u$.  Let $T$ be an $(n-1,n-1)$-current on
$D$.

The {\it slice} of $T$ with respect to a hyperplane $z_1=a$ is the
current
 $$
\langle T,z_1,a\rangle (\psi)
=\lim_{\epsilon\to 0}\frac{1}{\pi\epsilon^2}
\int_{\{|z_1-a|\le\epsilon\}\cap D}\psi(z_2,\ldots,z_n)\frac{1}{4}
dd^c|z_1|^2\wedge T.
 $$
Here $\psi$ is a $C_0^\infty$ test form on $z_1=a$, extended to $D$
independently of $z_1$.

Now let $u_1,\ldots,u_{n-1}$ and $w$ be bounded pluri\-sub\-har\-monic
functions on $D$, and put $T=w\,dd^cu_1\wedge\cdots\wedge
dd^cu_{n-1}$. Then by \cite{BT3, Proposition 4.1}, $\langle
T,z_1,a\rangle $ exists for every $a\in\Bbb C$ and
 $$
\langle T,z_1,a\rangle
=\frac{1}{2\pi}w(a,z')dd^c u_1(a,z')\wedge\cdots\wedge
dd^cu_{n-1}(a,z').
 $$
 Here $z'=(z_2,\ldots,z_n)$. 

Finally, if $F$ is holomorphic on $D$ and $M=\{z\in D:F(z)=0\}$, then
by changing variables and since only regular points of $M$ have to be
taken into account, one gets
$$\langle w(dd^cu)^{n-1}, F,0\rangle = w|_M(dd^cu)^{n-1}.$$

We write $\epsilon'=(\epsilon_2,\ldots,\epsilon_n)$, 
${\epsilon'}^2=\prod_{j=2}^n \epsilon^2_j$,
and $|z'|<\epsilon'$ for $|z_j|<\epsilon_j$, $j=2,\ldots,n$.

\proclaim{Lemma 3.5} Let $\psi=\psi(z_1)$ be a test function on
$\{z\in D: z'=0'\}$, and let $w$ and $u$ be bounded
pluri\-sub\-har\-monic functions on a bounded domain $D\subset\Bbb
C^n$. Then
 $$\split
&\int_{\{z_2=0,\ldots,z_n=0\}}\psi(z_1)w(z_1,0')\, dd^c u(z_1,0')\\
=
&\lim_{\epsilon'\downarrow 0}\frac{1}{2^{n-1}{\epsilon'}^2}
\int_{\{|z'|<\epsilon'\}}\psi(z_1)w(z)
dd^c|z_2|^2\wedge\ldots\wedge dd^c|z_n|^2 \wedge d d^cu.
\endsplit\tag3.5$$
\endproclaim

\demo{Proof} Apply slicing with respect to $z_2=0$ to the current
$T=w\,dd^c|z_3|^2\wedge\ldots\wedge dd^c|z_n|^2\wedge dd^cu$ to obtain
$\langle T, z_2,0\rangle=\frac{1}{2\pi}w(z_1,0,z_3,\ldots,
z_n)dd^c|z_3|^2\wedge\ldots\wedge dd^c|z_n|^2\wedge dd^cu$.  Next in
$\{z_2=0\}$, apply slicing with respect to $z_3=0$ to the current
$w\,dd^c|z_4|^2\wedge\ldots\wedge dd^c|z_n|^2\wedge dd^cu$.
Continuing in this fashion we obtain (3.5).
\enddemo

\demo{Proof of the `if part' of Theorem 3.1} We keep our notation and
proceed as follows. First we will show that $dd^c f\ge 0$ on the
compact neighborhood $K=K_{z_0}$ of $z_0$ provided by Lemma 3.3. Next
we apply Lemma 3.5 to show that the restriction of $f$ to any complex
line passing through $z_0$ is finely subharmonic on a fine
neighborhood of $z_0$.

Let $v$ be any plurisubharmonic function on a ball $B$ in $\Bbb C^n$,
let $h$ be a $\Bbb C$-affine bijection of $\Bbb C^n$, and let $\phi\in
C^\infty_0(B)$ be a test function. Then the action of the Riesz
measure $\Delta(v\circ h)$ on $\phi\circ h$ can be expressed as
follows
 $$
\split 4^{n-1}(n-1)!&\int_{h^{-1}(B)}\phi\circ h (z) \Delta(v\circ
h)\\
 &=\int_{h^{-1}(B)}\phi\circ h (z) dd^c(v\circ
h)\wedge\left(dd^c\|z\|^2\right)^{n-1}\\
 &=\int_B\phi(\zeta)dd^cv(\zeta)\wedge
\left(dd^c\|h^{-1}(\zeta)\|^2\right)^{n-1}.
 \endsplit$$

Returning to $f$, we have by Lemma 3.2 that the Riesz measure $\Delta
(f\circ h)$ is positive on $h^{-1}(K)$, hence (with $h^{-1}=g$) we
obtain that
 $$
dd^cf(\zeta)\wedge\left(dd^c\|g(\zeta)\|^2\right)^{n-1}\tag3.6
 $$
 is a positive measure on $K$ for every $\Bbb C$-affine bijection $g$
of $\Bbb C^n$, and by continuity also for every $\Bbb C$-affine map
$g:\Bbb C^n\to\Bbb C^n$.

To finish the proof we want to show that $f$ restricted to a complex
line $L$ passing through $z_0$ is finely subharmonic in a fine
neighborhood of $z_0$ relative to $L$.  We write $z=(z_1, z')$ and can
assume that $z_0=0$ and that $L$ is given by $z'=0'$.  Because $K$ is
an $\Cal F$-neighborhood of $0$, there exists a bounded non-negative
plurisubharmonic function $w$ defined on a ball $B_0$ about 0 that
equals 0 on $B_0\setminus K$, while $w(0)>0$. Then $\{z\in B_0:
w(z)>0\}$ is an $\Cal F$-open subset of $K$ that contains $0$. On $K$
we have $f=f_1-f_2$ where $f_1, f_2$ are plurisubharmonic on a ball
containing $K$, hence on $B_0$. We apply Lemma 3.5 to $f_1$ and $f_2$
separately and subtract to obtain from (3.6) (with
$g(z)=g(z_1,z')=(0,z')$)

 $$\split
&\int_{\{z_2=0,\ldots,z_n=0\}}\psi(z_1)w(z_1,0')\, dd^c f(z_1,0')\\
=&\lim_{\epsilon'\downarrow 0}\frac{1}{2^{n-1}{\epsilon'}^2}
\int_{\{|z'|<\epsilon'\}}\psi(z_1)w(z)
dd^c|z_2|^2\wedge\ldots\wedge dd^c|z_n|^2 \wedge d d^cf.
 \endsplit$$
If $\epsilon'$ is sufficiently small then the integrals occurring in the
limit on the right hand side are non-negative for every non-negative
test function $\psi$ on $B_0\cap \{z'=0'\}$.

We conclude that the Riesz measure of $f|_L$ is positive on a
neighborhood of $0$. By Lemma 3.2, $f|_L$ is finely subharmonic on
this neighborhood. Varying $z_0$ over $L$ and using the
sheaf-property, we find that $f|L$ is finely subharmonic.
\hfill{$\square$}\enddemo

\proclaim{Corollary 3.6} Let $f$ be a bounded weakly $\Cal
F$-pluri\-sub\-har\-monic function on an $\Cal F$-domain
$\Omega\subset\Bbb C^n$ such that $f$ admits the representation
$f=f_1-f_2$ of Lemma 3.3 on $\Omega$, and let $\chi_K$ denote the
characteristic function of a compact set $K$ in $\Omega$. Then for
$\Bbb C$-affine functions $l_1,\ldots, l_{n-1}$ on $K=K_{z_0}$ from
Lemma 3.3
 $$
\int_\Omega\chi_K(z) dd^c|l_1|^2\wedge \cdots
\wedge dd^c|l_{n-1}|^2\wedge dd^cf\ge 0.
 $$
\endproclaim
\demo{Proof}
This follows from (3.6) with
$g(z)=(l_1(z),\cdots,\l_{n-1}(z),0)$.
 \hfill{$\square$}
 \enddemo

From Theorem 3.1 we derive the following two results, one about
removable singularities for weakly $\Cal F$-pluri\-sub\-har\-monic
functions, and the other about the supremum of a family of such
functions.

\proclaim{Theorem 3.7} Let $f:\Omega\to[-\infty,+\infty[\,$ be $\Cal
F$-locally bounded from above on an $\Cal F$-open set
$\Omega\subset\Bbb C^n$, and let $E$ be an $\Cal F$-closed pluripolar
subset of \,$\Omega$. If $f$ is weakly $\Cal
F$-pluri\-sub\-har\-monic on $\Omega\setminus E$ then $f$ has a unique
extension to a weakly $\Cal F$-pluri\-sub\-har\-monic function on all
of \,$\Omega$, and this extension $f^*$ is given by
 $$
 f^*(z)=\Cal F\text{-}\underset{\Sb\zeta\to z\\
 \zeta\in\Omega\setminus E\endSb}\to{\lim\sup}\,f(\zeta),\qquad
 z\in\Omega.
 $$
 \endproclaim

\demo{Proof} The stated function $f^*$ (the $\Cal F$-upper
semicontinuous regularization of $f$) equals $f$ on the $\Cal F$-open
set $\Omega\setminus E$ because $f$ is $\Cal F$-upper
semi\-con\-tinu\-ous on $\Omega\setminus E$. Furthermore, $f^*$ is
$\Cal F$-upper semi\-con\-tinu\-ous and $<+\infty$ on all of $\Omega$
(finiteness because $f$ is $\Cal F$-locally bounded from above). By
the `only if part' of Theorem 3.1, for any $\Bbb C$-affine bijection
$h$ of $\Bbb C^n$, $f^*\circ h$ therefore is $\Bbb R^{2n}$-finely
hypo\-har\-monic on $h^{-1}(\Omega\setminus E)=h^{-1}(\Omega)\setminus
h^{-1}(E)$ and $\Cal F$-upper semi\-con\-tinu\-ous $<+\infty$ on
$h^{-1}(\Omega)$. In particular, $f^*\circ h<+\infty$ is $\Bbb
R^{2n}$-finely upper semi\-con\-tinu\-ous on $h^{-1}(\Omega)$. Because
$E$ is pluripolar so is $h^{-1}(E)$, which thus is $\Bbb
R^{2n}$-polar.  According to \cite{F1, Theorem 9.14}, $f\circ h$ is
therefore $\Bbb R^{2n}$-finely hypoharmonic on all of
$h^{-1}(\Omega)$, and so $f^*$ is indeed weakly $\Cal
F$-pluri\-sub\-har\-monic on $\Omega$, by the `if part' of Theorem
3.1.  Because the pluri\-polar set $E$ has empty $\Cal F$-interior,
$f^*$ is the only weakly $\Cal F$-pluri\-sub\-har\-monic and hence
$\Cal F$-continuous extension of $f$ to $\Omega$.
{\hfill$\square$}\enddemo

In view of Lemma 2.8 there is a similar result about removable
singularities for weakly $\Cal F$-holo\-morph\-ic functions:

\proclaim{Corollary 3.8} Let $h:\Omega\to\Bbb C$ be $\Cal F$-locally
bounded on $\Omega$ $(\Cal F$-open in $\Bbb C^n)$. If $h$ is weakly
$\Cal F$-holo\-morph\-ic on $\Omega\setminus E$ $(E$ $\Cal F$-closed
and pluri\-polar in $\Bbb C^n)$ then $h$ extends uniquely to a weakly
$\Cal F$-holo\-morph\-ic function $h^*:\Omega\to\Bbb C$, given by
 $$
h^*(z)
=\Cal F\text{-}\underset{\Sb\zeta\to z\\\zeta\notin E\endSb}
\to\lim\,h(\zeta),\qquad z\in\Omega.
 $$
 \endproclaim

\proclaim{Theorem 3.9} Let $\Omega$ denote an $\Cal F$-open subset of
$\Bbb C^n$.  For any $\Cal F$-locally upper bounded family of weakly
$\Cal F$-pluri\-sub\-har\-monic functions $f_\alpha$ on $\Omega$, the
least $\Cal F$-upper semi\-con\-tinu\-ous majorant $f^*$ of the
pointwise supremum $f=\sup_\alpha f_\alpha$ is likewise weakly $\Cal
F$-pluri\-sub\-har\-monic on $\Omega$, and
$\{z\in\Omega:f(z)<f^*(z)\}$ is pluripolar.
 \endproclaim

\demo{Proof} We may assume that the set $A$ of indices $\alpha$ is
upper directed and that the net $(f_\alpha)_{\alpha\in A}$ is
increasing; furthermore that $\Omega$ is $\Cal F$-connected and that
$f_\alpha\not\equiv-\infty$ for some $\alpha\in A$.  For any function
$f:\Omega\to[-\infty,+\infty[\,$ which is $\Cal F$-locally bounded
from above, write
 $$
f^*(z)=\Cal F{\text-}\underset{\zeta\to z}\to{\lim\sup}\,f(\zeta),
\quad\check f(z)=\Bbb R^{2n}{\text-}\underset{\zeta\to
z}\to{\fine\lim\sup}\,f(\zeta).
 $$
Then $\check f(z)\le f^*(z)<+\infty$, the former inequality because the
$\Bbb R^{2n}$-fine topology is finer than the $\Cal F$-topology.

As in Theorem 3.1, let $h:\Bbb C^n\to\Bbb C^n$ be a $\Bbb C$-affine
bijection, and note that
 $$ 
f\circ h=\sup_\alpha(f_\alpha\circ h),\quad(f\circ
h)\check{\vphantom{f}}=\check f\circ h,\qquad\text{on
}\;h^{-1}(\Omega),
 $$
the latter equation because $h$ is an $\Bbb R^{2n}$-fine
homeomorphism.  By Theorem 3.1, $f_\alpha\circ h$ is $\Bbb
R^{2n}$-finely hypoharmonic. Now $f_\alpha\circ h\le f^*\circ
h$. Furthermore, $f^*$ and hence $f^*\circ h$ and $\check f\circ h$
are $\Bbb R^{2n}$-finely locally bounded from above. It follows by
\cite{F1, Lemma 11.2} that $\check f\vphantom{\int^|}\circ h=(f\circ
h)\check{\vphantom{f}}\vphantom{\int^|}$ is $\Bbb R^{2n}$-finely
hypoharmonic.

We proceed to show that $\check f= f^*$ on $\Omega$, and hence that $\check
f$ is $\Cal F$-upper semicontinuous there. Invoking also Theorem 3.1
we shall thus altogether find that $\check f= f^*$ becomes $\Cal
F$-pluri\-sub\-har\-monic on $\Omega$, and in particular $\Cal
F$-continuous there, by Theorem 2.4(c).

Consider a point $z_0\in\Omega$ such that $f(z_0)>-\infty$. Fix
$\beta\in A$ with $f_\beta(z_0)>-\infty$, and choose an $\Cal F$-open
$\Cal F$-neighborhood $U$ of $z_0$ so that $U\subset\Omega$ and
 $$
f_\beta(z_0)-1<f_\beta\le f^*<f^*(z_0)+1\qquad\text{on }U,
 $$
noting that the weakly $\Cal F$-pluri\-sub\-har\-monic function
$f_\beta$ is $\Cal F$-continuous and that $f^*$ is $\Cal F$-upper
semicontinuous and $<+\infty$. Since $f_\beta\le f_\alpha\le f$ for
every $\alpha\curlyeqsucc\beta$ in $A$, any such $f_\alpha$ maps $U$
into some fixed bounded interval. According to Theorem 2.4(a),(b)
there exist $r>0$, an $\Cal F$-open set $O$ such that $z_0\in O\subset
B(z_0,r)$, and locally bounded ordinary pluri\-sub\-har\-monic
functions $\phi_\alpha$ and $\psi$ on $B(z_0,r)$ such that
$f_\alpha=\phi_\alpha-\psi$ on $O$ for every $\alpha\curlyeqsucc\beta$
in $A$.  The net $(\phi_\alpha)$ is increasing, along with the given
net $(f_\alpha)$. The pluri\-sub\-har\-monic functions $\phi_\alpha$
and $\psi$ are $\Cal F$-continuous, in particular $\Bbb R^{2n}$-finely
continuous. Writing $\sup_\alpha\phi_\alpha=\phi$ and denoting by
$\bar\phi$ the Euclidean $\Bbb R^{2n}$-subharmonic regularization of
$\phi$ in $B(z_0,r)$, we therefore have $\check\phi=\bar\phi$ there,
by Brelot's fundamental convergence theorem, see e.g.\ \cite{D,
1.XI.7}. Because $\check\phi\le\phi^*\le\bar\phi$ it follows that
$\check\phi=\phi^*$ in $B(z_0,r)$, and consequently
 $$
\check f=(\phi-\psi)\check{\vphantom|}=\check\phi-\psi=\phi^*-\psi
=(\phi-\psi)^*=f^*\quad\text{on }O
 $$
since $\psi$ is $\Cal F$-continuous and hence $\Bbb R^{2n}$-finely
continuous on the $\Cal F$-open, hence $\Bbb R^{2n}$-finely open set
$O\subset B(z_0,r)$.

Next, the set $\{z\in O:f(z)<f^*(z)\}=\{z\in O:\phi(z)<\phi^*(z)\}$ is
pluripolar, by the deep theorem of Bedford and Taylor \cite{BT1}, or
see \cite{K Theorem 4.7.6}. Writing 
 $$
E=\{z\in\Omega:f(z)<f^*(z)\},\quad e=\{z\in\Omega:f(z)=-\infty\},
 $$
we have thus found that every point $z_0\in\Omega\setminus e$ has an
$\Cal F$-neighborhood $O\subset\Omega\setminus e$ for which $O\cap E$
is pluripolar. Because
$e=\bigcap_{\alpha\in A}\{z\in\Omega:f_\alpha(z)=-\infty\}$ is $\Cal
F$-closed relative to $\Omega$, and pluripolar (some $f_\alpha$ being
$\not\equiv-\infty$), we infer by the quasi-Lindel\"of principle
\cite{BT2, Theorem 2.7} that indeed $E$ is pluripolar. Finally, we
have found that $f^*$ is $\Cal F$-pluri\-sub\-har\-monic on each $\Cal
F$-open set $O$ as above (as $z_0$ varies), and hence on their union
$\Omega\setminus e$, by the sheaf property. Because $f^*$ is $\Cal
F$-upper semicontinuous and $<+\infty$ on $\Omega$, and that $e$ is
pluripolar, we conclude from Theorem 3.7 above that indeed $f^*$ is
weakly $\Cal F$-pluri\-sub\-har\-monic on all of $\Omega$.
{\hfill$\square$}\enddemo

Taking for $\Omega$ a {\it Euclidean} open set we obtain in particular the
following

\proclaim{Corollary 3.10} For any family $\{f_\alpha\}$ of ordinary
pluri\-sub\-har\-monic functions on a Euclidean open set
$\Omega\subset\Bbb C^n$ such that $f:=\sup_\alpha f_\alpha$ is locally
bounded from above, the least pluri\-sub\-har\-monic majorant of $f$
exists and can be expressed as the upper semicontinuous regularization
of $f$ in the Euclidean topology on $\Bbb C^n$, as well as in the $\Cal
F$-topology and in the $\Bbb R^{2n}$-fine topology; that is, $\bar
f=f^*=\check f$.
 \endproclaim

The version of this involving the {\it Euclidean} topology is due to
Lelong \cite{L1}, or see \cite{L2, p.\ 26} or \cite{K, Theorem
2.9.10}. Being locally bounded from above, $f$ is in particular $\Cal
F$-locally bounded from above, and hence so is $f^*$, which is $\Cal
F$-pluri\-sub\-har\-monic by Theorem 3.9. Because $\Omega$ is
Euclidean open, it follows by Proposition 2.14 that $f^*$ even is an
ordinary pluri\-sub\-har\-monic function. From $f\le f^*\le\bar f$ it
therefore follows that $f^*=\bar f$. Similarly, $\check f=\bar f$ in
view of \cite{F1, Theorem 9.8(a)}.

The identity $f^*=\bar f$ is perhaps new even in the Euclidean case.

\smallskip 
We close this section with an alternative proof of the \lq only if
part' of Theorem 3.1. It is a bit shorter than the proof given
above. On the other hand it draws substantially on the theory of
functions of Beppo Levi and Deny, cf.~\cite{DL}, and its connection to
fine potential theory, cf.~\cite{F4}. We will need this approach again
in Section 4.

Following Deny \cite{DL} and subsequently \cite{F4} we consider for a
given Greenian domain $D$ (denoted $\Omega$ in \cite{DL} and
\cite{F4}) of $\Bbb C^n\cong\Bbb R^{2n}$ the complex Hilbert space
 $$
\wh{\Cal D}^1(D),
 $$
the completion of $\Cal D(D)=C^\infty_0(D,\Bbb C)$ in the Dirichlet
norm $||u||_1=||\nabla u||_{L^2(D,\Bbb C)}$. (For $n\ge2$ we may thus
take $D=\Bbb C^n$. For $n=1$, any bounded domain $D$ will do.)  Note
that $\wh{\Cal D}^1(D)$ is a space of distributions, \cite{DL,
Th\'eor\graveaccent eme 2.1 (p.\ 350)}. Elements of $\wh{\Cal D}^1(D)$
may be represented by quasi-continuous functions that are finite
quasi-everywhere.  For an $\Bbb R^2$-finely open set $\Omega\subset D$
denote by $\wh{\Cal D}^1(D, \Omega)$ the Hilbert subspace consisting
of all $\phi\in\wh{\Cal D}^1(D)$ such that some (and hence any) $\Bbb
R^{2n}$-quasi-continuous representative of $\phi$ satisfies $\phi=0$
$\Bbb R^{2n}$-quasi-everywhere on $D\setminus\Omega$, cf.\ \cite{DL,
Th\'eor\graveaccent eme 5.1, pp.\ 358\,f.}.  The positive cone in for
example $\wh{\Cal D}^1(D,\Omega)$ is denoted by $\wh{\Cal
D}^1_+(D,\Omega)$.  Let $V_l$ denote Lebesgue measure on $\Bbb C^l$,
and write $V_n=V$.

According to \cite{F4, Th\'eor\graveaccent eme 11} an $\Bbb
R^{2n}$-finely continuous (hence quasicontinuous) function
$f\in\wh{\Cal D}^1(D)$ is finely subharmonic quasi-everywhere (hence
actually everywhere by \cite{F1, Theorem 9.14}) on $\Omega$, if and
only if $f<+\infty$ and the inequality sign holds in (3.7):
 $$
 \frac14\int_D\nabla f\cdot\nabla\phi\,dV
=\sum_{j=1}^n\int_D(\partial_jf)\,(\bar\partial_j\phi)dV\le0\tag3.7
 $$
for every $\phi\in\wh{\Cal D}^1_+(D,\Omega)$. (It suffices of course to
integrate over $\Omega$.)

\demo{Alternative proof of the `only if part' of Theorem 3.1} Consider
a weakly $\Cal F$-pluri\-sub\-har\-monic function $f$ on an $\Cal
F$-open set $\Omega\subset\Bbb C^n$; hence $f$ is $\Cal F$-continuous
and $<+\infty$.  We leave out the trivial case $n=1$.  We may assume
that $f>-\infty$ on $\Omega$ (otherwise replace $f$ by $\max\{f,-p\}$
and let $p\to+\infty$).  It suffices to prove that $f$ is $\Bbb
R^{2n}$-finely hypo\-harmonic. 

Write $z=(z_1,\dots,z_n)=(z_1,z')\in\Bbb C^n$.  According to Theorem
2.4(a), every point $z_0\in\Omega$ then has an $\Cal F$-open $\Cal
F$-neighborhood $O \subset\Omega$ on which $f=f_1-f_2$, $f_1$ and
$f_2$ being bounded pluri\-sub\-har\-monic $>-\infty$ on some open
ball $B=B(z_0,r)$ containing $O $. In particular, $f_1$ and $f_2$ are
$\Bbb R^{2n}$-subharmonic on $B$. We may further assume that $-f_1$
and $-f_2$ are $\Bbb R^{2n}$-{\it potentials} on $B$, for otherwise we
may replace $-f_i$ for $i=1,2$ by its swept-out (relative to $B$) $\wh
R_{-f_i}^A$ on $A=B(z_0,r/2)$ (and $O $ by $O \cap A$). In terms of
the Green kernel $G$ on $B$ we may therefore write $-f_i=G\mu_i$ on
$B$ for some bounded positive measure $\mu_i$ of compact support in
$B$. Since $-f_i$ is bounded, its $G$-energy $\int G\mu_id\mu_i$ is
finite, and hence $G\mu_i$ is of Sobolev class
$W_0^{1,2}(B)\subset\wh{\Cal D}^1(\Bbb C^n,B)$, \cite{La, pp.\ 91--99},
cf.\ \cite{DL, Th\'eor\graveaccent eme 3.1 (p.\ 315)}.

For every $z'\in\Bbb C^{n-1}$ we have the $\Bbb C$-finely open set
 $$
O(z')=\{z_1\in\Bbb C:(z_1,z')\in O\}.
 $$
Because $f$ is weakly $\Cal F$-pluri\-sub\-har\-monic and $>-\infty$
on $O$, $f|_{L\cap O}$ is finely subharmonic for every complex line
$L$ in $\Bbb C^n$. It follows that (3.7) holds with $z$ replaced by
$z_1$ and with $O$ replaced by $O(z')$ for each $z'\in\Bbb C^{n-1}$:
 $$
 \int_{O(z')}\nabla_1 f(z_1, z')\cdot\nabla_1\phi(z_1,z')\,dV_1\le0.\tag3.8
 $$
Here $\nabla_1 =(\partial/\partial x_1,\partial/\partial y_1)$.
 Integrating (3.8)
with respect to $V_{n-1}$ leads by Fubini's theorem to
 $$
\int_O \nabla_1f(z_1,z')\cdot\nabla_1\phi(z_1,z')\,dV
\le0.
 $$
Similarly with the subscript $1$ replaced by any
$j\in\{1,\dots,n\}$. After addition this leads to
 $$
\int_O\nabla f\cdot\nabla\phi\,dV\le0.
 $$
According to \cite{F4, Th\'eor\graveaccent eme 11} quoted above,
this shows that $f$ indeed is $\Bbb R^{2n}$-finely \-sub\-har\-monic on
$O$, and hence, by varying $z_0$, on all of $\Omega$. 
 {\hfill$\square$}\enddemo

\heading{\bf 4. Biholomorphic invariance}\endheading

The sigma-algebra $\QB$ of {\it quasi Borel sets} in $\Bbb C^n$ is
generated by the Borel sets and the sets of capacity 0 (see
\cite{BT2}). $\QB$ contains the finely open sets. All currents
originating from wedge products of $dd^c$ of bounded plurisubharmonic
functions have measure coefficients that are Borel measures and put no
mass on pluripolar sets, hence they extend naturally to $\QB$.

\proclaim{Proposition 4.1} Let $f$ be a bounded weakly $\Cal
F$-pluri\-sub\-har\-monic function on an $\Cal F$-domain
$\Omega\subset\Bbb C^n$ such that $f$ admits the representation
$f=f_1-f_2$ of Lemma 3.3 on $\Omega$, and let $\chi_K$ denote the
characteristic function of a compact set $K\subset\Omega$. Then for
holomorphic functions $g_1,\ldots, g_{n-1}$ on $K=K_{z_0}$ from Lemma
3.3
 $$
\int_\Omega\chi_K(z) dd^c|g_1|^2\wedge \cdots
\wedge dd^c|g_{n-1}|^2\wedge dd^cf\ge 0.\tag4.1
 $$
\endproclaim

\demo{Proof} Corollary 3.6 yields that (4.1) is valid for
compact sets $\tilde K\subset K_{z_0}$ and $\Bbb C$-affine functions
$g_i$. For arbitrary holomorphic functions $g_j$ we have
 $$\split
&\int_\Omega\chi_K(z)dd^c|g_1|^2\wedge\cdots
\wedge dd^c|g_{n-1}|^2\wedge dd^cf\\
&=\lim_{N\to\infty}\sum_{j=1}^N
\int_\Omega\chi_{E_j^N}(z)dd^c|l^{j,N}_1|^2
\wedge\cdots\wedge dd^c|l^{j,N}_{n-1}|^2
\wedge dd^cf
\endsplit\tag4.2$$
for suitable quasi Borel sets $E_j^N$ and complex affine approximants
$l^{j,N}_k$ of $g_k$ on $E_j^N$ ($k=1,\ldots,n$).  Hence the right
hand side of (4.2) is indeed non-negative.  {\hfill$\square$}\enddemo

\proclaim{Theorem 4.2} Let $\Omega$ be $\Cal F$-open in $\Bbb C^n$
$(n\ge2)$. Given an $\Cal F$-continuous function $f\in \wh{\Cal
D}^1(\Bbb C^n)$ with values in $[-\infty,+\infty[$\,, the following
are equivalent:

{\rm (a)} $f$ is weakly $\Cal F$-pluri\-sub\-harmonic on $\Omega$,

{\rm (b)} for every $\phi\in\wh{\Cal D}^1_+(\Bbb C^n,\Omega)$ and every
$\lambda=(\lambda_1,\dots,\lambda_n)\in\Bbb C^n$,
$$
\sum_{j,k=1}^n \lambda_j\bar\lambda_k\int_\Omega
(\partial_jf)(\bar\partial_k\phi)\,dV\le0,
$$

{\rm (c)}  for every regular holomorphic map $h:\omega\to\Bbb C^n$ $(\omega$
open in $\Bbb C^n)$, $f\circ h$ is weakly $\Cal
F$-pluri\-sub\-har\-monic on $h^{-1}(\Omega)$ $(\subset\omega)$.
\endproclaim

\demo{Proof} (a)$\Rightarrow$(b). Using the characterization of weakly
$\Cal F$-pluri\-sub\-har\-monic functions given in Theorem 3.1, one
may adapt the proof of the `only if part' of \cite{K, Theorem 2.9.12}
as follows.  Suppose $f\in\wh{\Cal D}^1(\Bbb C^n)$ is weakly $\Cal
F$-pluri\-sub\-har\-monic on $\Omega$, and so $f\circ T$ is $\Bbb
R^{2n}$-finely sub\-har\-monic on $T^{-1}(\Omega)$ for any $\Bbb C$-affine
bi\-ject\-ion $T$ of $\Bbb C^n$. To prove (b) with constant
$\lambda=(\lambda_1,\ldots,\lambda_n)\in\Bbb C^n$, take
 $$
T_\epsilon(z)=z_1\lambda+\epsilon\sum_{l=2}^nz_le_l,\qquad\epsilon>0,
 $$
where $(e_1,\ldots,e_n)$ denotes the canonical base of $\Bbb
C^n$. From (3.7) we obtain (with integrations over $\Bbb C^n$),
replacing $\Omega$ and $\phi$, as we may, by $T_\epsilon^{-1}(\Omega)$
and $\phi\circ T_\epsilon\in\wh{\Cal D}^1_+(\Bbb
C^n,T_\epsilon^{-1}(\Omega))$,
 $$\align 0
&\ge\sum_{l=1}^n\int\partial_l(f\circ T_\epsilon)
 \bar\partial_l(\phi\circ T_\epsilon)\,dV\\
&=\sum_{j,k=1}^n\int[(\partial_jf)\circ T_\epsilon]
 [(\bar\partial_k\phi)\circ T_\epsilon](\lambda_j\bar\lambda_k
 +O(\epsilon))\,dV\\
&=|\det
 T_\epsilon|^2\biggl(\sum_{j,k=1}^n\int(\partial_jf)
 (\bar\partial_k\phi) \lambda_j\bar\lambda_k\,dV+O(\epsilon)\biggr).
 \endalign$$
 This leads to (b) after division by $|\det T_\epsilon|^2$ when we make
$\epsilon\to0$.

(b)$\Rightarrow$(a). Consider any  $\Bbb C$-affine bijection
$T=(T_1,\ldots,T_n)$ of $\Bbb C^n$, say
 $$
T_l(z)=\sum_{j=1}^nc_{lj}z_j+d_l,\qquad l\in\{1,\ldots,n\},\,z\in\Bbb C,
 $$ 
 with $c_{lj},d_l\in\Bbb C$ and $\det T\ne0$. We obtain
 $$\align
 \int\partial_l(f\circ T)\bar\partial_l(\phi\circ T)dV
&=\sum_{j,k=1}^n\int[(\partial_jf)\circ T][(\bar\partial_k\phi)\circ
 T] c_{lj}\bar c_{lk}\,dV\\
&=|\det T|^2\sum_{j,k=1}^n\int c_{lj}\bar
c_{lk}(\partial_jf)(\bar\partial_k\phi)\,dV\le0
 \endalign$$
by (b) with $\lambda_j=c_{lj}$.  After division by $|\det T|^2$ and
summation over $l$ this shows according to (3.7) and Theorem 3.1 that
the $\Cal F$-continuous function $f<+\infty$ indeed is $\Cal
F$-pluri\-sub\-har\-monic on $\Omega$.

(c)$\Rightarrow$(a).  This is contained in Theorem 3.1 (even with $h$
in (c) just a $\Bbb C$-affine bijection and with $f\circ h$ just $\Bbb
R^{2n}$-finely subharmonic).

(a)$\Rightarrow$(c). We may assume that $f>-\infty$ on $\Omega$
(otherwise pass to $f_p:=\max\{f,-p\}$, $p\in\Bbb N$, and let
$p\to+\infty$). According to Theorem 2.4(a), every point $z_0\in
h^{-1}(\Omega)$ then has an $\Cal F$-open $\Cal F$-neighborhood
$O\subset h^{-1}(\Omega)$ on which $f=f_1-f_2$, $f_1$ and $f_2$ being
bounded pluri\-sub\-har\-monic on some open set $D\subset\Bbb C^n$
containing $O$. In particular, $\Omega$ and $O$ are $\Bbb
R^{2n}$-finely open, and $f_1$ and $f_2$ are $\Bbb R^{2n}$-subharmonic
on $D$. We may further assume that the Jacobian matrix
$(\partial_jh_k)$ of the regular holomorphic map $h:\omega\to\Bbb C^n$
is bounded with determinant bounded away from $0$.

Denoting by $\Cal S(D,O)$ the convex cone of all functions of class
$\wh{\Cal D}^1(D)$ which are $\Bbb R^{2n}$-finely superharmonic
quasi-everywhere on $O$, we have by \cite{F4, p.\ 129} that $-f\in
\Cal S(D,O)$ and hence by \cite{F4, Th\'eor\graveaccent eme 11(b)}
 $$
\int_O\nabla f\cdot\nabla\phi\,\,dV\le0\qquad\text{for
}\phi\in\wh{\Cal D}_+^1(D,O).
 $$
For any $\psi\in\wh{\Cal D}^1(\Bbb C^n)$ we have (by the properties of
$h$ required above) $\psi\circ h\in\wh{\Cal D}^1(\omega)$. According
to Theorem 3.1 it suffices to show that the $\Cal F$-continuous
function $f\circ h$ is $\Bbb R^{2n}$-finely subharmonic on
$h^{-1}(O)$. For this it suffices by (3.7) to prove that, for every
$j\in\{1,\ldots,n\}$,
 $$
\int_{h^{-1}(O)}(\partial_j(f\circ
h))(\bar\partial_j\psi)\,dV\le0\qquad\text{for every }\psi\in
\wh{\Cal D}^1_+(\Bbb C^n,h^{-1}(O)),
 $$
and here $\psi$ may be replaced equivalently by $\phi\circ h$ with
$\phi\in\wh{\Cal D}^1_+(D,O)$ (or just as well with $\phi\in\Cal
D_+(D,O)$). We take $j=1$ and write
 $$
dV=(i/2)^n\,dz_1\wedge d\bar z_1
\wedge\ldots\wedge dz_n\wedge d\bar z_n=(1/4)^n
dd^c|z_1|^2\wedge\ldots\wedge dd^c|z_n|^2.
 $$
Then we obtain by the chain rule, writing by abuse of notation
$h^{-1}=(h_1^{-1},h_2^{-1},\allowmathbreak\ldots, h^{-1}_n)$ in terms
of the inverse $h^{-1}$ of the map $h$
 $$\allowdisplaybreaks\align
&\int_{h^{-1}(O)}(\partial_1(f\circ
h))(\bar\partial_1(\phi\circ h))dd^c|z_1|^2\wedge\ldots\wedge
dd^c|z_n|^2\\
=&\int d(f\circ h)\wedge d^c(\phi\circ h)\wedge
dd^c|z_2|^2 \wedge\ldots\wedge dd^c| z_n|^2\\ =&\int d(f\circ h)\wedge
d^c(\phi\circ h) \wedge dd^c|h^{-1}_2\circ h |^2 \wedge\ldots \wedge
dd^c|h^{-1}_n\circ h |^2 \\
=&\int df\wedge d^c\phi\wedge
dd^c|h^{-1}_2|^2 \wedge\ldots\wedge dd^c|h^{-1}_n|^2\tag4.3\\
=&\int\phi d(d^c f\wedge dd^c|h^{-1}_2|^2 \wedge\ldots\wedge
dd^c|h^{-1}_n|^2) \tag4.4\\
=&-\int\phi dd^c f\wedge
dd^c|h^{-1}_2|^2 \wedge\ldots\wedge dd^c|h^{-1}_n|^2.\\
 \endalign$$
The last three lines are in $h$-coordinates.  Equality (4.3) is
justified by approximating $f$ and $\phi$ in $\wh{\Cal D}^1$ with
functions in $\Cal D$ and applying Stokes' theorem to the approximants.
The final expression is non-positive because of Proposition 4.1, and
we are done. 
 {\hfill$\square$}\enddemo

Now we wish to consider the case where $h$ is just some sort of
pluri\-finely holo\-morph\-ic map. Recall from the text preceding
Proposition 2.9 that an $n$-tuple $(h_1,\ldots,h_n)$ of
strongly/weakly $\Cal F$-holo\-morph\-ic functions $h_j:U\to\Bbb C$
($U$ $\Cal F$-open in some $\Bbb C^m$) is termed a {\it
strongly/weakly} $\Cal F$-{\it holo\-morph\-ic map} (or {\it curve} if
$m=1$).

\definition{Definition 4.3} A {\it strongly $\Cal F$-biholomorphic
map} $h$ from an $\Cal F$-open set $U\subset\Bbb C^n$ onto its image in
$\Bbb C^n$ is an $\Cal F$-homeo\-morph\-ism with the property that there
exists for every $z\in U$ a compact $\Cal F$-neighborhood $K_z$ of $z$
in $U$ and a $C^\infty$-diffeo\-morph\-ism $\Phi_z$ from an open
neighborhood of $K_z$ to its image in $\Bbb C^n$ such that
$\Phi_z|_{K_z}=h|_{K_z}$ and that $\Phi_z|_{K_z}$ is a $C^1$-limit of
holo\-morphic maps defined on open sets containing $K_z$.
 \enddefinition

\proclaim{Proposition 4.4} The composition $f\circ h$ of a weakly $\Cal
F$-pluri\-sub\-har\-monic function $f$ on an $\Cal F$-open set
$\Omega\subset\Bbb C^n$ with a strongly $\Cal F$-bi\-holo\-morph\-ic
map $h:U\to\Omega $ $(U$ $\Cal F$-open in $\Bbb C^n)$ is weakly $\Cal
F$-pluri\-sub\-har\-monic on $h^{-1}(\Omega)$ $(\subset \Bbb C^n)$.
\endproclaim

\demo{Proof} For $n=1$ this is contained in \cite{F3, Corollaire, p.\
63} (in which $h$ is any finely holo\-morph\-ic function on
$U$). Suppose therefore that $n\ge2$. We may assume that $\Omega$ is
$\Cal F$-connected and that $f\not\equiv-\infty$, and so $f$ is in
particular $\Bbb R^{2n}$-finely subharmonic.  As shown in the
beginning of the alternative proof of Theorem 3.1 given at the end of
Section 3 we may further suppose that $\Omega$ is bounded in $\Bbb
C^n$ and that $f$ is bounded and of class $\wh{\Cal D}^1(D)$
for some bounded domain $D\subset\Bbb C^n$ containing $\Omega$.  Fix
$z\in U$ and let $K_z$ be a compact $\Cal F$-neighborhood of $z$ in
$U$ on which $h$ has the properties described in Definition 4.3.  It
will be sufficient to see that the expression (4.3) is non-positive if
$\phi\in\wh{\Cal D}_+^1(D,O)$ for some $\Cal F$-open set $O\subset D$
with $z\in O\subset K_z$.  Notice that $df\wedge d^c\phi$ is a form
with $L^1$ coefficients that is supported on $K_z$. Thus let $(h_m)$
be a sequence of bi-holomorphic maps on open sets containing $K_z$
that converge in $C^1$ to $h$ on $K_z$. Then
 $$ 
\lim_{m\to\infty}dd^c|h_{m,2}^{-1}|^2 \wedge\ldots\wedge
dd^c|h_{m,n}^{-1}|^2=dd^c|h^{-1}_2|^2 \wedge\ldots\wedge
dd^c|h^{-1}_n|^2,
 $$ 
uniformly on $K_z$. Now the expression (4.4) is non-positive when we
replace $h$ by $h_m$. By Lebesgue's dominated convergence theorem we
conclude that (4.4) is also non-positive for $h$ a strongly
$\Cal F$-biholomorphic map. 
 {\hfill$\square$}\enddemo

It is reasonable to expect that the concept of weakly $\Cal
F$-pluri\-sub\-har\-monic function is invariant even under composition with suitable
{\it weakly} $\Cal F$-biholomorphic mappings. Currently, we do not
know of a fine inverse function theorem for weakly $\Cal
F$-holomorphic maps of several variables. In fact we don't even know
if weakly $\Cal F$-holomorphic functions have weakly $\Cal
F$-holomorphic partial derivatives.  However, we can handle the
special case of a map of the form
 $$
 G(z)=g(z_1)+(0,z_2,\ldots,z_n),\tag4.5
 $$ 
where $g=(g_1,\dots,g_n)$ is a finely holomorphic curve in $\Bbb C^n$,
and that turns out to be sufficient.

\proclaim{Theorem 4.5} The composition $f\circ g$ of a weakly $\Cal
F$-plurisubharmonic $($resp.\ weakly $\Cal F$-holo\-morph\-ic$)$
function $f$ on an $\Cal F$-open set $\Omega\subset\Bbb C^n$ with a
finely holomorphic curve $g$ in $\Bbb C^n$ defined on a finely open
set $U\subset{\Bbb C}$, is finely hypoharmonic $($resp.\ finely
holo\-morph\-ic$)$ on the finely open pre-image $g^{-1}(\Omega)$
$(\subset\Bbb C)$.
\endproclaim

\demo{Proof} The theorem is known for $n=1$, cf.\ \cite{F3, \S4 and
Th\'eor\graveaccent eme 13(a)}, so we suppose that $n\ge2$. According
to Theorem 2.10(a), $g$ is continuous from $U$ with the fine topology
to $\Bbb C^n$ with the plurifine topology. The pre-image
$g^{-1}(\Omega)$ is therefore finely open in $\Bbb C$, and $f\circ g$
is finely upper semi\-con\-tinu\-ous (even finely continuous). We may
of course assume that $f$ is bounded, that $U$ is finely connected,
and that $g$ is non-constant, for example that $g_1$ is non-constant,
hence a fine-to-fine open map, \cite{F3, p.\ 64}.

Given a point $z_0\in g^{-1}(\Omega)$ ($\subset U\subset\Bbb C$) with
$g'_1(z_0)\ne0$, cf.\ \cite{F3, Corollaire 11}, there exists a finely
open set $O\subset\Bbb C$ such that $z_0\in O\subset g^{-1}(\Omega)$
and that $g_1$ is injective on $O$ with $g'_1(z)\ne0$ for every $z\in
O$, hence $(g_1|_O)^{-1}$ is finely holomorphic on the finely open set
$g_1(O)$, \cite{F3, Th\'eor\graveaccent eme 13}. We may further assume
after diminishing $O$ that there exists a $C^\infty$-map $\phi:\Bbb
C\to\Bbb C^n$ such that $\phi=g$ on $O$ and hence $\partial\phi=g'$,
$\bar\partial\phi=0$ on $O$, \cite{F3, Th\'eor\graveaccent eme
11(c)}. Since $\partial\phi_1(z_0)=g'_1(z_0)\ne0$ we may arrange (by
further diminishing $O$) that $\phi_1$ is injective on some open set
$\omega\subset\Bbb C$ containing the closure of $O$ in $\Bbb C$, and
hence that $\phi_1|_\omega$ is a $C^{\infty}$-diffeo\-morph\-ism of
$\omega$ onto $\phi_1(\omega)$. Likewise, we may achieve that there
exists a sequence of curves $\phi^{(\nu)}$ such that each coordinate
$\phi^{(\nu)}_j$ of $\phi^{(\nu)}$, $j\in\{1,\ldots,n\}$, is a
rational function defined on some open set $O^{(\nu)}\subset\omega$
containing $O$, and that $\phi^{\nu}\to\phi$,
$(\phi^{\nu})'\to\phi'$ ($=\partial\phi$), uniformly on $O$ as
$\nu\to\infty$, cf.\ \cite{F3, Theorem 11(a)}.  With this final choice
of $O$ define $G:O\times\Bbb C^{n-1}\to\Bbb C^n$ by (4.5), writing now
$t\in\Bbb C^n$ in place of $z\in\Bbb C^n$.

For a given point $z=(z_1,z')\in O\times\Bbb C^{n-1}$ choose a compact
fine neighborhood $L_{z_1}$ of $z_1$ in $\Bbb C$ so that
$L_{z_1}\subset O$ ($\subset\omega$), and a number
$c>\max\{|z_2|,\ldots,|z_n|\}$.  In analogy with (4.5) define
 $$
\Phi_z(t)=\phi(t_1)+(0,t_2,\ldots,t_n)\;
\qquad\text{for }t\in\omega\times\Bbb C^{n-1}.
 $$
It is easily verified that $\Phi_z$ is a $C^\infty$-diffeo\-morph\-ism
of $\omega\times\Bbb C^{n-1}$ onto its image in $\Bbb C^n$. We have
$\Phi_z=G$ on $O\times\Bbb C^{n-1}$, and in particular on the compact
$\Cal F$-neighborhood $K_z:=L_{z_1}\times[-c,c]^{n-1}$ of $z$ in
$O\times\Bbb C^{n-1}$, and consequently $G:O\times\Bbb C^{n-1}\to\Bbb
C^n$ is a strongly $\Cal F$-biholomorphic map.

Suppose first that $f$ is weakly {\it $\Cal F$-pluri\-sub\-har\-monic}
on $\Omega$.  According to Proposition 4.4, $f\circ G$ is a weakly
$\Cal F$-pluri\-sub\-har\-monic map defined on the $\Cal F$-open set
$G^{-1}(\Omega)\subset O\times\Bbb C^{n-1}$. Now
 $$
G(t_1,0,\ldots,0)=g(t_1)=(g_1(t_1),g_2(t_1),\ldots,g_n(t_1))
\qquad\text{for }t_1\in O, 
 $$
 and hence $(f\circ G)(t_1,0,\ldots,0)=(f\circ g)(t_1)$ for $t_1\in
O$. It follows that $f\circ g$ is finely hypoharmonic on $O$, hence so
(by varying $z_0$) on $\{t\in g^{-1}(\Omega):g'_1(t)\ne0\}$, which
differs only by a countable and hence polar set from $g^{-1}(\Omega)$,
cf.\ \cite{F3, Th\'eor\graveaccent eme 15}. Because $f\circ g$ is
finely continuous and $<+\infty$ on $g^{-1}(\Omega)$
($\subset U$) we conclude by the removable singularity theorem
\cite{F1, Theorem 9.15} that indeed $f\circ g$ is finely hypoharmonic
on $g^{-1}(\Omega)$.

Finally, let instead $f$ be weakly {\it $\Cal F$-holo\-morp\-hic} on
$\Omega$, in particular weakly (complex) $\Cal F$-har\-monic, by
\cite{F3, D\'efinition 3}.  As shown in the alternative proof of the
`only if part' of Theorem 3.1 given at the end of Section 3 we may
suppose that $\Omega$ is contained in some bounded domain
$D\subset\Bbb C^n$ and that $f$ is bounded and of class $\wh{\Cal
D}^1(D)$. Each of the functions $z\mapsto z_jf(z)$,
$j\in\{1,\ldots,n\}$, is therefore bounded and of class $\wh{\Cal
D}^1(D)$.

By the former part of the theorem, the bounded functions
$\pm\Re(f\circ g)$ and $\pm\Im(f\circ g)$ are finely subharmonic on
$g^{-1}(\Omega)$, and hence $f\circ g$ is (complex) finely harmonic
there. Similarly, $(z_j f)\circ g$ $=$ $g_j\cdot(f\circ g)$ is finely
harmonic on $g^{-1}(\Omega)$ ($\subset U\subset\Bbb C$),
$j\in\{1,\ldots,n\}$.  We may therefore further assume that $U$ is
contained in a bounded domain $D_1\subset\Bbb C$ and that $f\circ g$
and each $g_j$ are bounded and of class $\wh{\Cal
D}^1(D_1)$. Some component of $g$, say $g_1$, is
non-constant, and it therefore again follows from \cite{F3,
Th\'eor\`eme 13} that every point $z_0\in g^{-1}(\Omega)$ ($\subset
U$) with $g'_1(z_0)\ne0$ has a finely open fine neighborhood $O\subset
g^{-1}(\Omega)$ on which $g_1$ is injective with $g'_1\ne0$; and the
inverse $h_1:=(g_1|_O)^{-1}$ is finely holomorphic on the finely open
set $g_1(O)$. So is therefore the re-parametrized curve $g\circ
h_1:g_1(O)\to\Bbb C^n$. The function $h$ on $g_1(O)\times\Bbb C^{n-1}$
defined by
 $$
h(z)=h(z_1,\ldots,z_n)=h_1(z_1)
 $$
 is weakly $\Cal F$-holo\-morph\-ic. So is therefore $hf$, and
consequently $(hf)\circ g$ is finely harmonic on $O$, like
$f\circ g$ above.  Note that $h\circ g=h_1\circ g_1$. For $t\in O$,
 $$
[(hf)\circ g](t)=[(h\circ g)(t)][(f\circ g)(t)]
=[(h_1\circ g_1)(t)][(f\circ g)(t)]=t\cdot(f\circ g)(t),
 $$
and since this function $(hf)\circ g=t\cdot(f\circ g)$ of $t\in O$ is
finely harmonic it follows according to the result of Lyons \cite{Ly},
cf.\ \cite{F3, \S3}, which was utilized in Lemma 2.8, that
$f\circ g$ is finely holomorphic on $O$, hence (by varying $z_0$)
quasi-everywhere on $\{t\in g^{-1}(\Omega):g'_1(t)\ne0\}$, and indeed
everywhere on $g^{-1}(\Omega)$ by the removable singularity theorem
for finely holomorphic functions, \cite{F3, Corollaire 3}.
{\hfill$\square$}\enddemo

The following extension of Theorem 4.5 from finely holo\-morph\-ic
curves $g$ to $\Cal F$-holomorphic {\it maps} $h$ is a strengthening
of Theorem 2.10(b),(c), in which $f$ was required to be {\it strongly}
$\Cal F$-pluri\-sub\-har\-monic (resp.\ strongly $\Cal
F$-holo\-morph\-ic). Likewise, Theorem 4.6 (for a weakly $\Cal
F$-pluri\-sub\-har\-monic function $f$) extends Proposition 4.4 (in
which $m=n$, and $h$ is strongly bi\-holo\-morph\-ic).

\proclaim{Theorem 4.6} The composition $f\circ h$ of a weakly $\Cal
F$-pluri\-sub\-har\-monic $($resp.\ weakly $\Cal F$-holo\-morph\-ic$)$
function $f$ on an $\Cal F$-open set $\Omega\subset\Bbb C^n$ with a
weakly $\Cal F$-holo\-morphic map $h:U\to\Bbb C^n$ $(U$ $\Cal F$-open in
$\Bbb C^m)$ is weakly $\Cal F$-pluri\-sub\-har\-monic $($resp.\ weakly
$\Cal F$-holo\-morph\-ic$)$ on $h^{-1}(\Omega)$ $(\subset\Bbb C^m)$.
 \endproclaim

\demo{Proof} According to Theorem 2.10(a), $h$ is continuous from
$U\subset\Bbb C^m$ to $\Bbb C^n$ with their respective plurifine
topo\-logies. It follows that $h^{-1}(\Omega)$ is $\Cal F$-open in
$\Bbb C^m$, and that $f\circ h$ is $\Cal F$-upper
semi\-con\-tinu\-ous. Next, we restrict $f\circ h$ to a complex line
$L$ in $\Bbb C^m$, and observe that $h|_{L\cap U}$ is a finely
holomorphic curve. By Theorem 4.5, $f\circ h$ restricted to $L$
therefore is finely hypo\-har\-monic (resp.\ finely holo\-morph\-ic),
and we are done.  {\hfill$\square$}\enddemo

\heading{\bf References}\endheading

{\eightpoint

\widestnumber\key{ElKa\,}

\ref\key BH\by Bliedtner, J., and W. Hansen\paper Simplicial cones
in potential theory, II (approximation theorems)\jour Inventiones
Math.\vol 46\yr1978\pages255--275\endref

\ref\key BT1\manyby Bedford, E., and B.A. Taylor\paper A new capacity
for pluri\-sub\-har\-monic functions\jour Acta Math.\vol
149\yr1982\pages1--40\endref

\ref\key BT2\bysame\paper Fine topology, \v Silov boundary, and
$(dd^c)^n$\jour J. Functional Anal.\vol72\yr1987\pages225--251\endref

\ref\key BT3\bysame\paper Plurisubharmonic functions with logarithmic
singularities\jour Ann.  Inst. Fourier
(Grenoble)\vol38\issue4\yr1988\pages133--171\endref

\ref\key DL\by Deny, J., et J. L. Lions\paper Les espaces de type de
Beppo Levi\jour Ann. Inst. Fourier
(Grenoble)\vol5\yr1953-54\pages305--370\endref

\ref\key D\by Doob, J.L.\book Classical Potential Theory and Its
Probabilistic Counterpart\bookinfo Grundlehren\vol 262\publ
dddfSpringer\publaddr Berlin\yr 1984
\endref

\ref\key EEW\by Edigarian, A. El Marzguioui, S., and J.
Wiegerinck\paper The image of a finely holomorphic map is pluripolar,
\jour Ann. Polon. Math.\vol97\issue2\yr2010 \pages137--149\endref

\ref\key ElKa\by El Kadiri, M.\paper Fonctions finement
plurisousharmonique et topologie fine\jour Rend. Accad. Naz. Sci. XL
Mem. Mat. Appl.(5) \vol27 \yr2003\pages77-88
\endref

\ref\key EW1\manyby El Marzguioui, S., and J. Wiegerinck\paper The
pluri-fine topology is locally connected\jour Potential
Anal.\vol245\yr2006\pages283--288\endref

\ref\key EW2\bysame\paper Connectedness in the pluri\-fine
topology\paperinfo Functional Analysis and Complex Analysis
(Proceedings, 2008)\inbook Contemporary Math.\vol 481\publ
Amer. Math. Soc.\publaddr Providence, R.I., U.S.A.\yr2009\endref

\ref\key EW3\bysame\paper Continuity properties of finely
pluri\-sub\-har\-monic functions and pluripolarity\paperinfo To appear
in Indiana Univ. Math. J.  http://arxiv.org/abs/0906.2081\endref

\ref\key F1\manyby Fuglede, B.\book Finely Harmonic Functions\bookinfo
Lecture Notes in Math., No. 289\yr1972\publ Spring\-er\publaddr
Berlin\endref

\ref\key F2\bysame\paper Fonctions harmoniques et fonctions finement
harmoniques\jour Ann. Inst. Fourier
(Grenoble)\vol24\issue4\yr1974\pages77--91\endref

\ref\key F3\bysame\paper Sur les fonctions finement holomorphes\jour
Ann. Inst. Fourier (Grenoble)\vol31\issue4\yr1981\pages57--88\endref

\ref\key F4\bysame\paper Fonctions BLD et fonctions finement
surharmoniques\inbook S\'eminaire de Th\'eorie du Potentiel No. 6,
Lecture Notes in Mathematics No. 906\publ Springer\publaddr
Berlin\yr1982\endref

\ref\key F5\bysame\paper Localization in fine potential theory and
uniform approximation by sub\-harm\-onic functions\jour J. Functional
Anal.\vol49\yr1982\pages57--72\endref

\ref\key F6\bysame\paper Fonctions finement holomorphes de plusieurs
variables -- un essai\inbook S\'eminaire d'Analyse
P. Lelong--P. Dolbeault--H. Skoda, 1983/85, pp. 133--145, Lecture
Notes in Math. 1198, Springer, Berlin, 1986\endref

\ref\key F7\bysame\paper Finely ho\-lo\-morph\-ic functions. A
survey\jour Revue Roumaine Math. Pures
\vol33\yr1988\pages283--295\endref

\ref\key F8\bysame\paper Harmonic morphisms applied to classical
potential theory\jour Nagoya Math. J. (to appear)\endref

\ref\key J\by Josefson, B.\paper On the equivalence between locally
polar and globally polar sets for pluri\-sub\-har\-monic functions on
$C^n$\jour Arkiv Mat.\vol16\yr1978\pages 109--115\endref

\ref\key K\by Klimek, M.\book Pluripotential Theory\publ Clarendon
Press\publaddr Oxford\yr1991\endref

\ref\key La\by Landkof, N.S.\book Foundations of Modern Potential Theory
\bookinfo Grundlehren\vol180\publ Springer\publaddr Berlin\yr1972\endref

\ref\key Le1\manyby Lelong, P.\paper Les fonctions
pluri\-sous\-har\-moniques\jour Ann. Sci. \'Ecole
Norm. Sup.\vol62\yr1945\pages301--328\endref

\ref\key Le2\bysame\book Fonctions plurisousharmoniques et
formes diff\'erentielles positives\publ Gordon \& Breach\publaddr New
York\yr1968\endref

\ref\key Ly\by Lyons, T. J.\paper Finely holomorphic
functions\jour J. Functional Anal.\vol 37\yr1980\pages1--18\endref}

\medskip
{\smc Facult\'e des Sciences

Universit\'e Mohamed V

B.P. 726, Sal\'e-Tabriquet

11010 Sal\'e 

Morocco}

elkadiri\@fsr.ac.ma

\medskip
{\smc Institute for Mathematical Sciences

University of Copenhagen

Universitetsparken 5

2100 Copenhagen 

Denmark}

fuglede\@math.ku.dk

\medskip
{\smc KdV Institute for Mathematics 

University of Amsterdam

Science Park 904

P.O. box 94248, 1090 GE Amsterdam

The Netherlands}

j.j.o.o.wiegerinck\@uva.nl

\enddocument